\documentclass[11pt,twoside]{amsart}
\usepackage{amsmath, amsthm, amscd, amsfonts, amssymb, graphicx, color}
\usepackage[bookmarksnumbered, colorlinks, plainpages]{hyperref}

\newtheorem{thm}{Theorem}[section]
\newtheorem{cor}[thm]{Corollary}

\newtheorem{defn}[thm]{Definition}

\numberwithin{equation}{section}
\def\pn{\par\noindent}

\begin{document}

\leftline{ \scriptsize \it \bf}

\vspace{1.3 cm}

\title{Cauchy's Equations and Ulam's Problem}
\author{ali sadeghi}

\thanks{{\scriptsize
\hskip -0.4 true cm MSC(2010): Primary: 39B52, 39B72; Secondary:
47H09, 47H15.
\newline Keywords: Cauchy functional equations, Stability, Supertability, Ulam's Problem.\\
\newline\indent{\scriptsize}}}

\maketitle

\begin{abstract}
Our aim is to study the Ulam's problem for Cauchy's functional
equations. First, we present some new results about the
superstability and stability of Cauchy exponential functional
equation and its Pexiderized for class functions on commutative
semigroup to unitary complex Banach algebra. In connection with
the problem of Th. M. Rassias and our results, we generalize the
theorem of Baker and theorem of L. Sz$\acute{e}$kelyhidi. Then
the superstability of Cauchy additive functional equation can be
prove for complex valued functions on commutative semigroup under
some suitable conditions. This result is applied to the study of
a superstability result for the logarithmic functional equation,
and to give a partial affirmative answer to problem $18$, in the
thirty-first ISFE. The hyperstability and asymptotic behaviors of
Cauchy additive functional equation and its Pexiderized can be
study for functions on commutative semigroup to a complex normed
linear space under some suitable conditions. As some consequences
of our results, we give some generalizations of Skof's theorem,
S.-M. Joung's theorem, and another affirmative answer to problem
$18$, in the thirty-first ISFE. Also we study the stability of
Cauchy linear equation in general form and in connection with the
problem of G. L. Forti, in the 13th ICFEI (2009), we consider
some systems of homogeneous linear equations and our aim is to
establish some common Hyers-Ulam-Rassias stability for these
systems of functional equations and presenting some applications
of these results.
\end{abstract}

\vskip 0.2 true cm


\pagestyle{myheadings}
\markboth{\rightline {\scriptsize  Ali Sadeghi}}
         {\leftline{\scriptsize Cauchy's Equations and Ulam's Problem}}

\bigskip
\bigskip


\section{\bf Introduction}
\subsection{Cauchy's Functional Equations and Ulam's Problem: Stability of Functional Equations}
In general, a functional equation is any equation that specifies
a function or some functions in implicit form, where an implicit
equation is a relation of the form $F(x_{1},..., x_{n}) = 0$,
where $F$ is a function of several variables (often a
polynomial). However, in this note, we deal with Cauchy's
functional equations in one and two variables that specifies at
most three functions. According to the our study and for the
readers convenience and explicit later use, we give the following
definitions for functional equations and some types of stability.
In our study, we will provide some reasons for these definitions.
Note that in this note, all theorms, where has previously been
proven by some authors, we will present all this theorems based of
our notions and definitions.
\begin{defn}
Let $S$ and $B$ be nonempty sets. Let $g_{i}:S^{2}\rightarrow S$
for $i\in \{1,...,4\}$ and $G_{j}:B^{2} \rightarrow B$ for $j\in
\{1,2\}$ be functions. Then the following equation
\begin{equation}\label{0.1}
  G_{1}[f(g_{1}(x,y)),f(g_{2}(x,y))]=G_{2}[f(g_{3}(x,y)),f(g_{4}(x,y))]
\end{equation}
is called functional equation that specifies unknown functions
$f:S\rightarrow B$. We use the symbol $\Im(S,B)$ for equation
(\ref{0.1}) as a functional equation that specifies unknown
functions $f:S\rightarrow B$ and denote the set all solutions of
this functional equation with $Z_{\Im(S,B)}$. Similarly, we can
give to definition of functional equation in single variable.
\end{defn}
\begin{defn}\label{hu-stable}
Let $S$ and $B$ be nonempty sets such that $(B,d)$ be a metric
space. Assuming that $\varphi:S^{2}\rightarrow [0,\infty)$ and
$\phi:S\rightarrow [0,\infty)$ are functions and $\Im(S,B)$ is a
functional equation (of the form equation (\ref{0.1})). If for
every functions $f:S\rightarrow B$ satisfying the inequality
\begin{equation}\label{0.2}
  d(G_{1}[f(g_{1}(x,y)),f(g_{2}(x,y))],G_{2}[f(g_{1}(x,y)),f(g_{2}(x,y))])\leq \varphi(x,y)
\end{equation}
for all $x, y\in S$, there exists $T\in Z_{\Im(S,B)}$ such that
\begin{equation}\label{0.3}
  d(f(x),T(x))\leq \phi(x)
\end{equation}
for all $x\in S$, then we say that the functional equation
$\Im(S,B)$ is Hyers-Ulam-Rassias stable on $(S,B)$ with control
functions $(\varphi, \phi)$ and we denoted it by "HUR-stable" on
$(S,B)$ with controls $(\varphi, \phi)$. Also we call the
function $T$ as "HUR-stable function". If $\varphi(x,y)$ in
(\ref{0.2}) and $\phi(x)$ in (\ref{0.3}) are replaced by real's
$\delta>0$ and $\varepsilon>0$ respectively, then we say that
corresponding phenomenon of the functional equation $\Im(S,B)$ is
the Hyers-Ulam stable on $(S,B)$ and we denoted it by "HU-stable"
on $(S,B)$ with controls $(\delta, \varepsilon)$, and we call the
function $T$ as "HU-stable function". Similarly, we can define of
above concepts for functional equations in single variable.
\end{defn}
\begin{defn}\label{superstable}
Let $S$ and $B$ be nonempty sets such that $(B,d)$ be a metric
space. Assuming that $\varphi:S^{2}\rightarrow [0,\infty)$ and
$\phi:S\rightarrow [0,\infty)$ are functions and $\Im(S,B)$ is a
functional equation. If for every functions $f:S\rightarrow B$
satisfying the inequality
\begin{equation}\label{0.4}
  d(G_{1}[f(g_{1}(x,y)),f(g_{2}(x,y))],G_{2}[f(g_{3}(x,y)),f(g_{4}(x,y))])\leq \varphi(x,y)
\end{equation}
for all $x, y\in S$, then either $f\in Z_{\Im(S,B)}$ or
\begin{equation}\label{0.5}
  d(f(x))\leq \phi(x)
\end{equation}
for all $x\in S$, then we say that the functional equation
$\Im(S,B)$ is "superstable" on $(S,B)$ with control functions
$(\varphi, \phi)$.
\end{defn}
\begin{defn}
Let $S$ and $B$ be nonempty sets such that $(B,d)$ be a metric
space. Assuming that $\varphi:S^{2}\rightarrow [0,\infty)$ is a
function and $\Im(S,B)$ is a functional equation. If for every
functions $f:S\rightarrow B$ satisfying the inequality
\begin{equation}\label{0.4}
  d(G_{1}[f(g_{1}(x,y)),f(g_{2}(x,y))],G_{2}[f(g_{3}(x,y))],f(g_{4}(x,y))])\leq \varphi(x,y)
\end{equation}
for all $x, y\in S$, then $f\in Z_{\Im(S,B)}$, then we say that
the functional equation $\Im(S,B)$ is "hyperstable" on $(S,B)$
with control $(\varphi)$.
\end{defn}
In early $19$th century, Cauchy has determined the general
continuous solution of each of the functional equations
\begin{equation}\label{1}
  f(x+y)=f(x)+f(y),
\end{equation}
\begin{equation}\label{2}
  f(x+y)=f(x)\cdot f(y),
\end{equation}
\begin{equation}\label{3}
  f(x\cdot y)=f(x)+f(y),
\end{equation}
\begin{equation}\label{4}
  f(x\cdot y)=f(x)\cdot f(y),
\end{equation}
for real-valued functions on some subsets of $\mathbb{R}$ (is the
real numbers field) and showing that they are, respectivel, $ax$,
$a^{x}$, $a\log x$, $x^{a}$, where in each case $a$ is an
arbitrary constant. These functional equations, its often called
the Cauchy additive, exponential, logarithmic, and multiplicative
functional equations respectively, in honor of A. L. Cauchy. For
some first generalization of these functional equations see
\cite{Carmichael} and for recent achievements see
\cite{Kannappan}. If no further conditions are imposed on $f$,
then (assuming the axiom of choice) there are infinitely many
other functions that satisfy these functional equations. This was
proved in 1905 by Georg Hamel using Hamel bases. Such functions
are sometimes called Hamel functions. Fore more information see
(\cite{Kuczma}, Chap. 5). The fifth problem on Hilbert's list is
a generalisation of equation (\ref{1}). Functions where there
exists a real number $c$ such that $f(cx)\neq cf(x)$ are known as
Cauchy-Hamel functions and are used in Dehn-Hadwiger invariants
which are used in the extension of Hilbert's third problem from
$3$-D to higher dimensions \cite{Boltianskii}, also see
\cite{Sah}. The properties of Cauchy's functional equations are
frequently applied to the development of theories of other
functional equations. Moreover, the properties of Cauchy's
functional equations are powerful tools in almost every field of
natural and social sciences.

The starting point of the stability theory of functional
equations was the problem formulated by S. M. Ulam in 1940 (see
\cite{Ulam}). More precisely, S. M. Ulam gave a wide-ranging talk
before a Mathematical Colloquium at the University of Wisconsin
in which he discussed a number of important unsolved problems.
Among those was the following question concerning the stability
of homomorphisms:

\emph{Let $(S,.)$ be a group and $(B,.,d)$ be a metric group, let
$\Im(S,B)$ be the functional equation $f(xy)=f(x)f(y)$. Does for
every $\varepsilon>0$, there exists a $\delta>0$ such that the
functional equation $\Im(S,B)$ is HU-stable on $(S,B)$ with
controls $(\delta, \varepsilon)$.}

For the first in 1941, D. H.
Hyers \cite{Hyers} gave an affirmative partial answer to this
problem for the case where $S$ and $B$ are assumed to be Banach
spaces. The result of Hyers is stated as follows:
\begin{thm}\label{Hyers}
Suppose that $S$ and $B$ are two real Banach spaces and $\Im(S,B)$
be the functional equation $f(x+y)=f(x)+f(y)$. Then for every
$\epsilon>0$, the functional equation $\Im(S,B)$ is HU-stable on
$(S,B)$ with controls $(\epsilon, \epsilon)$ and also proved that
HU-stable function is unique.
\end{thm}
This is the reason for which today this type of stability is
called Hyers-Ulam stability of functional equations. And also this
is a reason for definition of (\ref{hu-stable}). Here, note that
Hyers considered only bounded control functions for Cauchy
additive functional equation. T. Aoki \cite{Aoki} introduced
unbounded one and generalized a result of Theorem ({\ref{Hyers}).
Th. M. Rassias, who independently introduced the unbounded
control functions was the first to prove the stability of the
linear mapping between Banach spaces. Taking this fact into
account, the additive functional equation $f(x+y)=f(x)+f(y)$ is
said to have the Hyers-Ulam-Rassias stability on two Banach
spaces, and this is a reason for definition of (\ref{hu-stable}).
This terminology is also applied to the case of other functional
equations. For more detailed definitions of such terminology one
can refer to \cite{Forti} and \cite{Hyers4}. Th. M. Rassias
\cite{Rassias} generalized Hyers's Theorem as follows:
\begin{thm}\label{Rassias}
Suppose that $S$ and $B$ are two real Banach spaces and $\Im(S,B)$
be the functional equation $f(x+y)=f(x)+f(y)$. Then for every
$\epsilon>0$ and $0\leq p<1$, the functional equation $\Im(S,B)$
is HUR-stable on $(S,B)$ with controls
$(\epsilon(\|x\|^{p}+\|y\|^{p}),
\frac{2\epsilon}{2-2^{p}}\|x\|^{p})$ and also proved that
HUR-stable function is unique.
\end{thm}
Th. M. Rassias \cite{Rassias4} noticed that the proof of this
theorem also works for $p<0$. In fourth section, we show that for
this case the Cauchy additive functional equation is hyperstable
and also give some generalizations. Following Th. M. Rassias and
P. $\check{S}$emrl \cite{Rassias5} generalized the result of
(\ref{Rassias}) and obtained stability result for the case
$p\geq0$ and $p\neq 1$. For the case $p=1$, Z. Gajda in his paper
\cite{Gajda} showed that the theorem of Th. Rassias
(\ref{Rassias}) is false for some special control function and
give the following counterexample.
\begin{thm}\label{Gajda}
Let $S=B$ be real field $\mathbb{R}$ and
$\Im(\mathbb{R},\mathbb{R})$ be the functional equation
$f(x+y)=f(x)+f(y)$. Then for every $\theta>0$ there is no
constant $\delta\in [0,\infty)$ such that the functional equation
$\Im(\mathbb{R},\mathbb{R})$ is HUR-stable on
$(\mathbb{R},\mathbb{R})$ with controls $(\theta(|x|+|y|)),
\delta|x|)$.
\end{thm}
M. S. Moslehian and Th. M. Rassias \cite{Moslehian} generalized
the Theorem (\ref{Hyers}) and Theorem (\ref{Rassias}) in
non-Archimedean spaces. Also, H. G. Dales and M. S. Moslehian in
\cite{Moslehian1} introduced multi-normed spaces and study some
properties of multi-bounded mappings on such spaces. Then they
proved some generalized Hyers-Ulam-Rassias stability theorems
associated to the Cauchy additive functional equation for mappings
from linear spaces into multi-normed spaces. The Hyers-Ulam
stability of mappings is in development and several authors have
remarked interesting applications of this theory to various
mathematical problems. For the first, L. C$\check{a}$dariu and V.
Radu proved the Hyers-Ulam-Rassias stability of the additive
Cauchy equation by using the fixed point method (see
\cite{Cadariu} and \cite{Radu}). This method appears to be
powerful and successful. In fact the Hyers-Ulam stability has
been mainly used to study problems concerning approximate
isometries or quasi-isometries, the stability of Lorentz and
conformal mappings, the stability of stationary points, the
stability of convex mappings, or of homogeneous mappings, etc.
For more information about theory of stability of functional
equations one can refer to \cite{Jung1}, \cite{Moszner}, and
\cite{Rassias1}. For important specials functional equations, we
can refer to
\cite{Baak,Baker,Hyers1,Hyers3,Isac,Jarosz,Johnson,Rassias2}, and
\cite{Tabor}. The stability problem of functional equations has
been extended in various directions and study by several
mathematicians. So, we are necessary to introduce exact
definitions of some stability which is applicable to all
functional equations in this note at the first. So, we present
these notions and definitions based of paper \cite{Rassias1} and
book \cite{Jung1}.
\subsection{Superstability of Functional Equations: Cauchy Exponential Functional Equations }
In 1979, another type of stability was observed by J. Baker, J.
Lawrence and F. Zorzitto \cite{Baker2}. Indeed, they proved that
if a real-valued function $f$ on a real vector space $V$ satisfies
the functional inequality
\begin{equation}\label{zori}
  |f(x+y)-f(x)f(y)|\leq \epsilon
\end{equation}
for some $\delta>0$ and for all $x, y\in V$, then $f$ is either
bounded or exponential. In fact, they prove that Cauchy
exponential functional equation is superstable on
$(V,\mathbb{R})$. This result was the first result concerning the
superstability phenomenon of functional equations and the
definition of (\ref{superstable}) is based on this result. Later,
J. Baker \cite{Baker1} (see also \cite{Czerwik} and \cite{Ger})
generalized this famous result as follows:
\begin{thm}\label{b1}
Let $(S,\cdot)$ be semigroup, $B$ be the field of complex numbers,
and $\Im(S,\mathbb{C})$ be the functional equation $f(x\cdot
y)=f(x)f(y)$. Then for every $\varepsilon>0$, the functional
equation $\Im(S,\mathbb{C})$ is superstable on $(S,\mathbb{C})$
with controls $(\varepsilon,(1+\frac{\sqrt{1+4\varepsilon}}{2}))$
\end{thm}
In the proof of the preceding theorem, the multiplicative
property of the norm was crucial. Indeed, the proof above works
also for functions $f:S\rightarrow A$, where $A$ is a normed
algebra in which the norm is multiplicative, i.e.,
$\|xy\|=\|x\|\|y\|$ for all $x, y\in A$. Examples of such real
normed algebras are the quaternions and the Cayley numbers. In
the same paper Baker gives the following example to show that
this result fails if the algebra does not have the multiplicative
norm property. Let $\epsilon>0$, choose $\delta>0$ so that
$|\delta-\delta^{2}|=\epsilon$ and let $f:\mathbb{C}\rightarrow
\mathbb{C}\oplus \mathbb{C}$ be defined as
$$f(\lambda)=(e^{\lambda},\delta),\ \ \ \ \lambda\in C.$$
Then, with the nonmultiplicative norm given by
$\|(\lambda,\mu)\|=\max\{|\lambda|,|\mu|\}$, we have
$\|f(\lambda+\mu)-f(\lambda)f(\mu)\|=\epsilon$ for all complex
$\lambda$ and $\mu$, f is unbounded, but it is not true that
$f(\lambda+\mu)=f(\lambda)f(\mu)$ for all complex $\lambda$ and
$\mu$.

The result of Baker, Lawrence and Zorzitto \cite{Baker2} was
generalized by L. Sz$\acute{e}$kelyhidi \cite{Szekelyhidi} in
another way and he obtained the following result.
\begin{thm}\label{sze}
Let $(S,\cdot)$ be an Abelian group with identity $1$, $B$ be the
field of complex numbers, and let $\Im(S,\mathbb{C})$ be the
functional equation $f(1)f(x \cdot y)=g(x)f(y)$, where
$g:S\rightarrow \mathbb{C}$ is a function. Assume that $M_{1},
M_{2}: S\rightarrow [0,\infty)$ are two functions. Then there
exists $\delta>0$ such that the functional equation
$\Im(S,\mathbb{C})$ is superstable on $(S,\mathbb{C})$ with
controls $(\min\{M_{1}(x),M_{2}(y)\}, \delta)$ such that
$f(x)=g(x)f(1)$ for all $x\in S$.
\end{thm}
During the thirty-first International Symposium on Functional
Equations, Th. M. Rassias \cite{Rassias3} introduced the term
$\emph{mixed stability}$ of the function
$f:E\rightarrow\mathbb{R}$ (or $\mathbb{C}$), where $E$ is a
Banach space, with respect to two operations 'addition' and
'multiplication' among any two elements of the set
$\{x,y,f(x),f(y)\}$. Especially, he raised an open problem
concerning the behavior of solutions of the inequality
$$\|f(x.y)-f(x)f(y)\|\leq\theta(\|x\|^{p}+\|y\|^{p}).$$

In connection with this open problem, P. G$\breve{a}$vruta
\cite{Gavruta} gave an answer to the problem suggested by Rassias
concerning the mixed stability.
\begin{thm}\label{gavr}
Let $S$ and $B$ be a real normed space and a normed algebra with
multiplicative norm, respectively. Let $\Im(S,B)$ be the
functional equation $f(x+y)=f(x)f(y)$. Then for every $\theta>0$
and $p>0$, the functional equation $\Im(S,B)$ is superstable on
$(S,B)$ with controls $(\theta(\|x\|^{p}+\|y\|^{p}), g(x))$, where
$g(x)=\frac{1}{2}(2^{p}+\sqrt{4^{p}+8\theta})\|x\|^{p}$ with
$\|x\|\geq1$.
\end{thm}
In the second section, first we give another proof of Theorem
(\ref{b1}), where its important idea for other results. We study
the superstability of Cauchy exponential functional equation. As a
consequence of our results and in connection with problem of T.
H. Rassias, we extend the results of Baker and
Sz$\acute{e}$kelyhidi in unitary complex Banach algebra. Also we
present this result for the Pexiderized Cauchy exponential
equation. More precisely, we proved the superstability and
stability of Cauchy exponential functional equation and its
Pexiderized when the controls functions is not bounded.
Furthermore, we consider the superstability and stability for the
funcional equation of the form $f(x+y)=g(x)f(y)$, in which $f$ is
a function from a commutative semigroup to an complex Banach
space and $g$ is function from a commutative to complex field and
next we consider the superstability and stability for the
equations of the forms $f(x+y)=g(x)f(y)$ and $f(x+y)=g(x)h(x)$
when $f, g$ and $h$ are three functions from a commutative
semigroup to an unitary complex Banach algebra. Also this Results
is applied to the study of homogeneous functional equation and its
Pexiderized.

Also during the thirty-first International Symposium on Functional
Equations (ISFE), the following question arises. Let $(S,\cdot)$
be an arbitrary semigroup or group and let a mapping
$f:S\rightarrow \mathbb{R}$ (the set of reals) be such that the
set $\{f(x\cdot y)-f(x)-f(y)\ |\ x, y\in S\}$ is bounded. Is it
true that there is a mapping $T:S\rightarrow R$ that satisfies
\begin{eqnarray}\label{5}
  T(x\cdot y)-T (x)-T (y) &=& 0
  \end{eqnarray}
for all $x, y\in S$ and that the set $\{T(x)-f(x)\ |\ x\in S\}$ is
bounded? G. L. Forti in \cite{Forti1} gave a negative answer to
this problem. In third section the superstability of the Cauchy
equation (in the sense of additive) can be proved for complex
valued functions on commutative semigroup under some suitable
conditions and the result is applied to the study of a
superstability result for the logarithmic functional equation.
Furthermore, these results is partial affirmative answers to
problem $18$, in the thirty-first ISFE. Also in fourth section, we
give another partial affirmative answers to this problem under
some suitable conditions.
\subsection{Asymptotic Behavior of Functional Equations}
Several authors have used asymptotic conditions in stating
approximations to Cauchy's functional equation. P. D. T. A.
Elliott \cite{Elliott} showed that if the real function $f$
belongs to the class $L^{p}(0, z)$ for every $z\geq0$, where
$p\geq1$, and satisfies the asymptotic condition
$$\lim_{z\rightarrow \infty} \frac{\int_{0}^{z}\int_{0}^{z}|f(x+y)-f(x)-f(y)|^{p}dxdy}{z}=0,$$
then there is a constant $c$ such that $f(x)=cx$ almost
everywhere on $\mathbb{R}^{+}$. One of the theorems of J. R.
Alexander, C. E. Blair and L. A. Rubel \cite{Alexander} states
that if $f\in L^{1}(0,b)$ for all $b>0$, and if for almost all
$x>0$
$$\lim_{u\rightarrow \infty} \frac{\int_{0}^{y}[f(x+y)-f(x)-f(y)]dy}{u}=0,$$
then for some real number $c$, $f(x)=cx$ for almost all $x\geq0$.

F. Skof \cite{Skof} proved the following theorem and applied the
result to the study of an asymptotic behavior of additive
functions.
\begin{thm}
Let $E_{1}$ and $E_{2}$ be a normed space and a Banach space,
respectively. Given $a>0$, suppose a function $f:E_{1}\rightarrow
E_{2}$ satisfies the inequality
$$\|f(x+y)-f(x)-f(y)\|\leq \delta$$
for some $\delta>0$ and for all $x, y\in E_{1}$ with
$\|x\|+\|y\|>a$. Then there exists a unique additive function
$A:E_{1}\rightarrow E_{2}$ such that
$$\|f(x)-A(x)\|\leq 9\delta$$
for all $x\in E_{1}$.
\end{thm}
Using this theorem, F. Skof \cite{Skof} has studied an interesting
asymptotic behavior of additive functions as we see in the
following theorem.
\begin{thm}\label{11}
Let $E_{1}$ and $E_{2}$ be a normed space and a Banach space,
respectively. For a function $f:E_{1}\rightarrow E_{2}$ the
following two conditions are equivalent:
\begin{enumerate}
  \item $\|f(x+y)-f(x)-f(y)\|\rightarrow 0$ as $\|x\|+\|y\|\rightarrow \infty$;
  \item $f(x+y)-f(x)-f(y)=0$
\end{enumerate}
for all $x, y\in E_{1}$.
\end{thm}

S.-M. Joung \cite{Joung}, proved that the Hyers-Ulam stability for
Jensen's equation on a restricted domain and the result applied to
the study of an interesting asymptotic behavior of the additive
mappings. More precisely, he proved that a mapping
$f:E_{1}\rightarrow E_{2}$ satisfying $f(0)=0$ is additive if and
only if
\begin{enumerate}
  \item $\|2f(\frac{x+y}{2})-f(x)-f(y)\|\rightarrow 0$ as $\|x\|+\|y\|\rightarrow \infty$.
\end{enumerate}

In fourth section, we study the hyperstability of Cauchy additive
functional equation and its Pexiderized for class functions on
commutative semigroups to complex normed spaces. As a consequence
of our results, we peresent some new results about the asymptotic
behavior of Cauchy additive functional equation and its
Pexiderized. Also, we give a simple proofs of Skof's theorem
(\ref{11}) and S.-M. Joung's theorem and show that these results
is true when $E_{2}$ be a complex normed linear space.
Furthermore, we present some generalization of Skof's and S.-M.
Joung's theorems and give another affirmative answer to problem
$18$, in the thirty-first ISFE.
\subsection{Stability and Common Stability of Functional Equations}

Of the most importance is the linear functional equation or
Cauchy linear equation in general form
\begin{equation}\label{e1}
    f(\rho(x))=p(x)f(x)+q(x)
\end{equation}
where $\rho$, $p$ and $q$ are given functions on an interval $I$
and $f$ is unknown. When $q(x)\equiv 0$ this equation, i.e.,
\begin{equation}\label{e2}
    f(\rho(x))=p(x)f(x)
\end{equation}
is called homogeneous linear equation. We refer the reader to
\cite{Kuczma1} and \cite{Agarwal} for numerous results and
references concerning this equation and  its stability in the
sense of Ulam.

In $1991$ Baker \cite{Baker3} discussed Hyers-Ulam stability for
linear equations (\ref{e1}). More concretely, the Hyers-Ulam
stability and the generalized Hyers-Ulam-Rassias stability for
equation
\begin{equation}\label{e3}
    f(x+p)=kf(x)
\end{equation}
were discussed by Lee and Jun \cite{Lee}. Also the gamma
functional equation is a special form of homogeneous linear
equation (\ref{e2}) were discussed by S. M. Jung
\cite{Jung2,Jung3,Jung4} proved the modified Hyers-Ulam stability
of the gamma functional equation. Thereafter, the stability
problem of gamma functional equations has been extended and
studied by several mathematicians
\cite{Alzer,Barnes,Jun,Kim,Kim1}.

Assume that $S$ is a nonempty set, $F=\mathbb{Q}$, $\mathbb{R}$
or $\mathbb{C}$, $B$ is a Banach space over $F$,
$\psi:S\rightarrow \mathbb{R}^{+}$, $f,g:S\rightarrow B$,
$p:S\rightarrow K\backslash\{0\}$, $q:S\rightarrow B$ are
functions, and $\sigma:S\rightarrow S$ is a arbitrary map.

In the fifth section of this note, we present some results about
Hyers-Ulam-Rassias stability via a fixed point approach for the
linear functional equation in general form (\ref{e1}) and its
Pexiderized
\begin{equation}\label{e4}
    f(\rho(x))=p(x)g(x)+q(x)
\end{equation}
under some suitable conditions. Note that the main results of this
section can be applied to the well known stability results for the
gamma, beta, Abel, Schr$\ddot{o}$der, iterative and G-function
type's equations, and also to certain other forms.

During the $13$st International Conference on Functional Equations
and Inequalities (ICFEI) $2009$ , G. L. Forti posed following
problem (see [\cite{ICFEI}, pp. 144]).

\emph{ Consider functional equations of the form
\begin{equation}\label{d1}
    \sum_{i=1}^{n}a_{i}f(\sum_{k=1}^{n_{i}}b_{ik}x_{k})=0 \ \ \ \ \
    \sum_{i=1}^{n}a_{i}\neq 0
\end{equation}
and
\begin{equation}\label{d2}
    \sum_{i=1}^{m}\alpha_{i}f(\sum_{k=1}^{m_{i}}\beta_{ik}x_{k})=0 \ \ \ \ \
    \sum_{i=1}^{m}\beta_{i}\neq 0
\end{equation}
where all parameters are real's number and
$f:\mathbb{R}\rightarrow \mathbb{R}$. Assume that the two
functional equations are equivalent, i.e., they have the same set
of solutions. Can we say something about the common stability?
More precisely, if (\ref{d1}) is stable, what can we say about the
stability of (\ref{d2}). Under which additional conditions the
stability of (\ref{d1}) implies that of (\ref{d2})?}

In connection the above problem, we consider the term of common
stability for systems of functional equations. We give the
definition of this type stability in fifth section.
In connection with the problem of G. L. Forti, we consider some
systems of homogeneous linear equations and our aim is to
establish some common Hyers-Ulam-Rassias stability for these
systems of functional equations. As a consequence of these
results, we give some superstability results for the exponential
functional equation.

For the reader's convenience and explicit later use, we will
recall a fundamental results in fixed point theory.
\begin{defn}
The pair $(X, d)$ is called a generalized complete metric space if
$X$ is a nonempty set and $d:X^{2}\rightarrow [0,\infty]$
satisfies the following conditions:
\begin{enumerate}
  \item $d(x,y)\geq0$ and the equality holds if and only if $x=y$;
  \item $d(x,y)=d(y,x)$;
  \item $d(x,z)\leq d(x,y)+d(y,z)$;
  \item every d-Cauchy sequence in X is d-convergent.
\end{enumerate}
for all $x, y\in X$.
\end{defn}
Note that the distance between two points in a generalized metric
space is permitted to be infinity.
\begin{thm}\cite{Diaz}\label{3.1}
Let $(X, d)$ be a generalized complete metric space and
$J:X\rightarrow X$ be strictly contractive mapping with the
Lipschitz constant $L$. Then for each given element $x\in X$,
either
$$d(J^{n}(x),J^{n+1}(x))=\infty$$
for all nonnegative integers $n$ or there exists a positive
integer $n_{0}$ such that
\begin{enumerate}
  \item $d(J^{n}(x),J^{n+1}(x))<\infty$, for all $n\geq n_{0}$;
  \item the sequence $\{J^{n}(x)\}$ converges to a fixed point $y^{\ast}$ of $J$;
  \item $y^{\ast}$ is the unique fixed point of $J$ in the set $Y=\{y\in X\ :\ d(J^{n_{0}}(x),y)<\infty\}$;
  \item $d(y,y^{\ast})\leq\frac{1}{1-L}d(J(y),y)$.
\end{enumerate}
\end{thm}
\begin{thm}(Banach's contraction principle)\label{3.1.1}
Let (X, d) be a complete metric space and let $J:X\rightarrow X$
be strictly contractive mapping. Then
\begin{enumerate}
  \item the mapping J has a unique fixed point $x^{\ast}= J(x^{\ast})$;
  \item the fixed point $x^{\ast}$ is globally attractive, i.e.,
  $$\lim_{n\rightarrow \infty} J^{n}(x)=x^{\ast}$$
  for any starting point $x\in X$;
  \item one has the following estimation inequalities:
  $$d(J^{n}(x),x^{\ast})\leq L^{n}d(x,x^{\ast}),$$
  $$d(J^{n}(x),x^{\ast})\leq\frac{1}{1-L}d(J^{n}(x),J^{n+1}(x)),$$
  $$d(x,x^{\ast})\leq\frac{1}{1-L}d(J(x),x)$$
\end{enumerate}
for all nonnegative integers $n$ and all $x\in X$.
\end{thm}
\section{\bf Superstability of Cauchy Exponential Functional Equation}
Some of the results of this section was proved by Mohsen
Alimohammady and Ali Sadeghi, which was published see
\cite{sadeghi1}. In the first, we give another proof of Baker's
Theorem (\ref{b1}). In general, this method appears to be
powerful and successful for our aims of this section.

\textbf{\emph{Proof}.} Let $f:S\rightarrow \mathbb{C}$ be a
unbounded function. Assume that $a\in S$ such that $|f(a)|>1$. Let
$\Im^{'}(S,B)$ be the functional equation $f(a\cdot y)=f(a)f(y)$.

\textbf{Step 1.} We prove that the functional equation
$\Im^{'}(S,B)$ is HU-stable on $(S,B)$ with controls
$(\varepsilon,\frac{\varepsilon}{|f(a)|-1})$ and HU-stable
function is unique.

we have
\begin{equation}\label{b3}
     |f(a\cdot y)-f(a)f(y)|\leq\varepsilon
\end{equation}
for all $y\in S$. Let us consider the set $A:=\{u:S\rightarrow
\mathbb{C}\}$ and introduce a metric on $A$ as follows:
$$d(u,h)=\sup_{y\in S}\frac{|u(y)-h(y)|}{\varepsilon}.$$
It is easy to show that $(A, d)$ is a complete metric space. Now
we define the function $J:A\rightarrow A$ with
$$J(h(y))=\frac{1}{f(a)}h(a\cdot y)$$ for all $h\in A$ and $y\in S$. So
  \begin{eqnarray*}
    d(J(u),J(h)) &=& \sup_{y\in S}\frac{|u(a\cdot y)-h(a\cdot y)|}{|f(a)|\varepsilon }\\
      &\leq& L\sup_{y\in S}\frac{|u(a\cdot y)-h(a\cdot y)|}{\varepsilon}\leq L d(u,h)
  \end{eqnarray*}
for all $u, h\in A$, that is $J$ is a strictly contractive
selfmapping of $A$, with the Lipschitz constant
$L=\frac{1}{|f(a)|}$. From (\ref{b3}), we get
$$|\frac{f(a\cdot y)}{f(a)}-f(y)|\leq L\varepsilon,$$ for all $y\in S$, which says that $d(J(f),f)\leq
L<\infty$ . So by Theorem \ref{3.1.1}, there exists a mapping
$T:S\rightarrow \mathbb{C}$ such that
  \begin{enumerate}
  \item $T$ is a fixed point of $J$, i.e.,
  \begin{eqnarray}
    T(a\cdot y) &=& f(a)T(y)
  \end{eqnarray}
  for all $y\in S$. The mapping $T$ is a unique fixed point of
  $J$. So, $T\in \Im^{'}(S,B)$ and is unique HU-stable function.
  \item $d(J^{n}(f),T)\rightarrow 0$ as $n\rightarrow \infty$. This implies that
  $$T(y)=\lim_{n\rightarrow\infty}\frac{f(a^{n}\cdot y)}{f(a)^{n}}$$
  for all $y\in S$.
  \item $d(f,T)\leq\frac{1}{1-L}d(J(f),f)$, which implies that
  $d(f,T)\leq\frac{L}{1-L}$ or
\begin{equation}\label{b4}
     |f(y)-T(y)|\leq \frac{\varepsilon}{|f(a)|-1}.
\end{equation}
for all $y\in S$. The proof of step 1 is complete.
 \end{enumerate}

\textbf{Step 2.} We prove that the HU-stable function $T$ have
the property $T(x\cdot y)=T(x)f(y)$ and then $f(x\cdot
y)=f(x)f(y)$ for all $x,y\in S$.

Let $x,y\in S$ be two arbitrary fixed elements, we have
$$|f(a^{n}\cdot x\cdot y)-f(a^{n}\cdot x)f(y)|\leq\varepsilon$$ and dividing by $|f(a)|^{n}$,
$$|\frac{f(a^{n}\cdot x\cdot y)}{f(a)^{n}}-\frac{f(a^{n}\cdot x)}{f(a)^{n}}f(y)|\leq\frac{\varepsilon}{|f(a)|^{n}}$$ and
letting $n$ to infinity, we get
\begin{equation}
    T(x\cdot y)=T(x)f(y)
\end{equation} for all $x,y\in S$. Let $x,y,z\in S$ be arbitrary
elements, then
$$T(x\cdot y\cdot z)=T(x)f(y\cdot z)$$ and $$T(x\cdot y\cdot z)=T(x\cdot y)f(z)=T(x)f(y)f(z)$$ or
$$T(x)(f(y\cdot z)-f(y)f(z))=0$$ for all $x,y,z\in S$. Since $f$ is an
unbounded function, then from (\ref{b4}) implies that the function
$T$ is nonzero. Therefore, we have $f(x\cdot y)=f(x)f(y)$ for all
$x,y\in S$ (i.e., $f\in \Im(S,B)$). The proof is complete.\\

In the following, first, we study the stability and superstability
of the a Pexider type of Cauchy exponential functional equation
$$f(x\cdot y)=g(x)f(y),$$
for class functions $f$ on commutative semigroup to complex Banach
space and given complex-valued function $g$. For the readers
convenience and explicit later use in this section, we present
the some notions.
\begin{defn}\label{tnv}
Let $(S,\cdot)$ be a semigroup and let $g:S\rightarrow
\mathbb{C}$ and $\psi: S^{2}\rightarrow [0,\infty)$ be functions,
then we denote the set $N_{g,\psi}$ as the all $a\in S$, which
$|g(a)|>1$ and
\begin{equation}\label{33.1}
  \psi(x,y\cdot a) \leq\psi(x,y)
\end{equation}
for all $x,y\in S$.
\end{defn}
\begin{defn}\label{3.2}
Let $(S,\cdot)$ be a semigroup, let $B$ be a complex Banach
algebra with unit $1_{B}$, and let $f:S\rightarrow B$ be a
function, then denote the set $M_{f}$ as follows:
$$M_{f}=\{a\in S\ :\ f(a)\in \mathbb{C} \times\{1_{B}\}\}.$$
Also, we introduce the function $\widehat{f}:S\rightarrow
\mathbb{C}$, where $\widehat{f}(a)\times 1_{B}=f(a)$ if $a\in
M_{f}$ and $\widehat{f}(a)=1$ for another elemnts.
\end{defn}
\begin{thm}\label{3.3}
Let $(S,\cdot)$ be commutative semigroup and $B$ be a complex
Banach space, and let $\Im(S,B)$ be the functional equation
$f(x\cdot y)=g(x)f(y)$, where $g:S\rightarrow \mathbb{C}$ is a
given function. If $N_{g,\psi}\neq {\O}$, then the functional
equation $\Im(S,B)$ is HUR-stable on $(S,B)$ with controls
$$(\psi(x,y),\inf_{a\in N_{g,\psi}}\frac{\psi(a,y)}{|g(a)|-1})$$
such that the HUR-stable function such $T$ is unique and
\begin{equation}\label{33.2}
(g(x\cdot y)-g(x)g(y))T(z)=0
\end{equation}
for all $x,y, z\in S$.
\end{thm}
\textbf{\emph{Proof}.} Let $a\in N_{g,\psi}$ be fixed and let
$\Im^{a}(S,B)$ be the functional equation $f(a\cdot y)=g(a)f(y)$.

\textbf{Step (1).} We prove that the functional equation
$\Im^{a}(S,B)$ is HUR-stable on $(S,B)$ with controls
$(\psi(a,y),\frac{\psi(a,y)}{|g(a)|-1})$ and HUR-stable function
such $T_{a}$ is unique.

We have
\begin{eqnarray}\label{3.5}
  \|f(a\cdot y)-g(a)f(y)\| &\leq& \psi(a,y)
\end{eqnarray}
for all $y\in S$. Let us consider the set $A:=\{g:S\rightarrow
B\}$ and introduce the generalized metric on $A$:
$$d(g,h)=\sup_{\{y\in S\ ;\ \psi(a,y)\neq0\}}\frac{\|g(y)-h(y)\|}{\psi(a,y)}.$$
It is easy to show that $(A, d)$ is complete metric space. Now we
define the function $J_{a}:A\rightarrow A$ as follows:
$$J_{a}(h(y))=\frac{1}{g(a)}h(y\cdot a)$$
for all $h\in A$ and $y\in S$. So
  \begin{eqnarray*}
    d(J_{a}(u),J_{a}(h)) &=& \sup_{\{y\in S\ ;\ \psi(a,y)\neq0\}}\frac{\|u(y\cdot a)-h(y\cdot a)\|}{|g(a)|\psi(a,y)} \\
      &\leq& \sup_{\{y\in S\ ;\ \psi(a,y)\neq0\}}\frac{\|u(y\cdot a)-h(y\cdot a)\|}{|g(a)|\psi(a,y\cdot a)}\leq \frac{1}{|g(a)|}d(u,h)
  \end{eqnarray*}
for all $u, h\in A$, that is $J$ is a strictly contractive
selfmapping of $A$, with the Lipschitz constant
$\frac{1}{|g(a)|}$. From (\ref{3.5}), we get
$$\|\frac{f(y\cdot a)}{g(a)}-f(y)\|\leq \frac{\psi(a,y)}{|g(a)|}$$

for all $y\in S$, which says that $d(J(f),f)<
\frac{1}{|g(a)|}<\infty$. By Theorem (\ref{3.1}), there exists a
mapping $T_{a}:S\rightarrow B$ such that
  \begin{enumerate}
  \item $T_{a}$ is a fixed point of $J$, i.e.,
  \begin{eqnarray}
    T_{a}(y\cdot a) &=& g(a)T_{a}(y)
  \end{eqnarray}
  for all $y\in S$. The mapping $T_{a}$ is a unique fixed point of $J$ in the set $\tilde{A}=\{h\in A\ :\
  d(f,h)<\infty\}$. Hence, $T_{a}\in \Im^{a}(S,B)$ and is unique
  HUR-stable function.
  \item $d(J^{n}(f),T_{a})\rightarrow 0$ as $n\rightarrow \infty$. This implies that
  $$T_{a}(y)=\lim_{n\rightarrow\infty}\frac{f(y\cdot n^{a})}{g(a)^{n}}$$
  for all $x\in S$.
  \item $d(f,T_{a})\leq \frac{1}{1-\frac{1}{|g(a)|}}d(J(f),f)$, which implies,
  $$\|f(y)-T_{a}(y)\| \leq \frac{\psi(a,y)}{|g(a)|-1}$$
  for all $y\in S$ and the proof of step (1) is complete.
 \end{enumerate}

\textbf{Step (2).}  We prove that the functional equation
$\Im(S,B)$ is HUR-stable on $(S,B)$ with controls
$(\psi(x,y),\inf_{a\in N_{g,\psi}}\frac{\psi(a,y)}{|g(a)|-1})$
such that the HUR-stable function such $T$ is unique and
$(g(x\cdot y)-g(x)g(y))T(z)=0$.

From (\ref{3.5}), its easy to show that following inequality
\begin{eqnarray}\label{3.6}
  \|f(y\cdot a^{n})-g(a)^{n}f(y)\| &\leq& \sum_{i=0}^{n-1}\psi(a,y\cdot a^{i})|g(a)|^{n-1-i}
\end{eqnarray}
for all $y, a\in S$ and $n\in \mathbb{N}$. Now since
$$\psi(x,y\cdot a) \leq \psi(x,y)$$ for all $y\in S$
and all $a\in N_{g,\psi}$, so

\begin{equation}\label{3.666}
  \psi(a,y\cdot a^{m})\leq \psi(a,y)
\end{equation}
for all $x\in S$ and $m\in \mathbb{N}$, thus from (\ref{3.6}), we
obtain
\begin{eqnarray}\label{3.7}
  \|f(y\cdot a^{n})-g(a)^{n}f(y)\| &\leq& \psi(a,y)\frac{|g(a)|^{n}-1}{|g(a)|-1}
\end{eqnarray}
for all $y\in S$. Our aim is to prove that $T_{a}=T_{b}$ for each
$a, b\in N_{g,\psi}$. We have from inequality (\ref{3.7})
\begin{eqnarray}
  \|f(y\cdot a^{n})-g(a)^{n}f(y)\| &\leq& \psi(a,y\cdot a^{n})\frac{|g(a)|^{n}-1}{|g(a)|-1} \label{3.8} \\
  \|f(y\cdot b^{n})-g(b)^{n}f(y)\| &\leq&\psi(b,y\cdot b^{n})\frac{|g(b)|^{n}-1}{|g(b)|-1} \label{3.9}
\end{eqnarray}
for all $y\in S$. On the replacing $y$ by $y\cdot b^{n}$ in
(\ref{3.8}) and $y$ by $y\cdot a^{n}$ in (\ref{3.9}) and so from
(\ref{3.666}), we get
$$\|f(y\cdot (a\cdot b)^{n})-g(a)^{n}f(y\cdot b^{n})\|\leq \psi(a,y) \frac{|g(a)|^{n}-1}{|g(a)|-1}$$
$$\|f(y\cdot (a\cdot b)^{n})-g(b)^{n}f(y\cdot a^{n})\|\leq\psi(b,y) \frac{|g(b)|^{n}-1}{|g(b)|-1}.$$ Thus,
$$\|g(a)^{n}f(y\cdot b^{n})-g(b)^{n}f(y\cdot a^{n})\|\leq\psi(a,y)\frac{|g(a)|^{n}}{|g(a)|-1}+
\psi(b,y)\frac{|g(b)|^{n}}{|g(b)|-1}$$ and dividing by
$|g(a)^{n}g(b)^{n}|$
$$\|\frac{f(y\cdot a^{n})}{g(a)^{n}}-\frac{f(y\cdot b^{n})}{g(b)^{n}}\|\leq \frac{\psi(a,y)}{|g(b)|^{n}(|g(a)|-1)}+\frac{\psi(b,y)}{|g(a)|^{n}(|g(b)|-1)}$$
and since $|g(a)|>1$ for any $a\in N_{g,\psi}$, letting $n$ to
infinity, we obtain $T_{a}(y)=T_{b}(y)$ for all $y\in S$. So,
there a unique function $T$ such that $T = T_{a}$ for every $a\in
N_{g,\psi}$ and
$$\|f(y)-T(y)\|\leq\frac{\psi(a,y)}{|g(a)|-1}$$ for
all $y\in S$ and $a\in N_{g}$. Since $a\in N_{g,\psi}$ is an
arbitrary element, so
$$\|f(y)-T(y)\|\leq\inf_{a\in N_{g}} \frac{\psi(a,y)}{|g(a)|-1} $$
for all $y\in S$.

Let $x, y\in S$ and $a\in N_{g,\psi}$ be three arbitrary fixed
elements, we have
$$\|f(x\cdot y \cdot a^{n})-g(x)f(y \cdot a^{n})\|\leq\psi(x,y \cdot a^{n}),$$  from (\ref{3.666}) and dividing last inequality
with $|g(a)|^{n}$, we get
$$\|\frac{f(x\cdot y \cdot a^{n})}{g(a)^{n}}-g(y)\frac{f(x\cdot a^{n})}{g(a)^{n}}\|\leq\frac{\psi(x,y)}{|g(a)|^{n}}$$
and letting $n$ to infinity, we get $T(x\cdot y)=g(x)T(y)$.

Let $x, y, z\in S$ be arbitrary elements, then
$$T(x\cdot y \cdot z)=g(x\cdot y)T(z)$$
and
$$T(x\cdot y \cdot z)=g(x)T(y\cdot z)=g(x)g(y)T(z)$$
or $$(g(x\cdot y)-g(x)g(y))T(z)=0$$ for all $x, y, z\in S$. The
proof is complete.\\

In connection with the problem of Th. M. Rassias and Theorem
(\ref{3.3}), in the following, we prove some extensions of Baker's
theorem (\ref{b1}) and also we prove a generalized version of L.
Sz$\acute{e}$kelyhidi's theorem (\ref{sze}). Note that in the
Definition (\ref{tnv}), if the function $\psi$ is constant such
$\epsilon$, then either $f$ is bounded or unbounded function if
and only if either $f(N_{f,\epsilon})$ is bounded or unbounded
set.
\begin{cor}
Let $(S,\cdot)$ be commutative semigroup and let
$\Im(S,\mathbb{C})$ be the functional equation $f(x\cdot
y)=f(x)f(y)$. If $f:S\rightarrow \mathbb{C}$ is a function such
that $f(N_{f,\psi})$ is an unbounded set, then $f\in
Z_{\Im(S,\mathbb{C})}$.
\end{cor}
\textbf{\emph{Proof.}} In Theorem (\ref{3.3}), if we put $
B:=\mathbb{C}$ and $g:=f$, then we have
$$|f(y)-T(y)|\leq\inf_{a\in N_{f,\psi}}
\frac{\psi(a,y)}{|f(a)|-1}$$ and $$T(x\cdot y)=f(x)T(y)$$ for all
$x, y\in S$. Now if $f(N_{f,\psi})$ is an unbounded set, then
$f=T$ and so $f\in Z_{\Im(S,B)}$.
\begin{cor}
Let $(S,\cdot)$ be commutative semigroup and let
$\Im(S,\mathbb{C})$ be the functional equation $f(x\cdot
y)=f(x)f(y)$. If $f:S\rightarrow \mathbb{C}$ is a function and
there exist $a\in N_{f,\psi}$ such that $\psi(a,y)=0$ for any
$y\in S$ or there is $x_{0}\in S$ such that the following limit
exists and
$$\lim_{n\rightarrow \infty} \frac{f(x_{0}\cdot
a^{n})}{(f(a))^{n}}\neq0,$$ then $f\in Z_{\Im(S,\mathbb{C})}$.
\end{cor}
\textbf{\emph{Proof.}} In Theorem (\ref{3.3}), if we put $
B:=\mathbb{C}$ and $g(x):=f(x)$, then we have
$$|f(y)-T(y)|\leq\inf_{a\in N_{f,\psi}}
\frac{\psi(a,y)}{|f(a)|-1}$$ and $(f(x\cdot y)-f(x)f(y))T(z)=0$
for all $x, y, z\in S$. If $\psi(a,y)=0$ for any $y\in S$, then
$f=T$ and the proof is done. For another case, according to the
proof of Theorem (\ref{3.3}) and assumption, we have
$T(x_{0})\neq 0$ and so that $f(x\cdot y)=f(x)f(y)$ for all $x,
y\in S$, where that is $f\in Z_{\Im(S,\mathbb{C})}$.
\begin{cor}\label{homo3}
Let $(S,\cdot)$ be commutative semigroup and $B$ be a complex
Banach space, and let $\Im(S,B)$ be the functional equation
$f(x\cdot y)=g(x)f(y)$, where $g:S\rightarrow \mathbb{C}$ is a
given function. If $g(N_{g,\psi})$ is an unbounded set, then the
functional equation $\Im(S,\mathbb{C})$ is hyperstable on $(S,B)$
with control $(\psi(x,y))$ and $f(x)=g(x)f(1)$ for all $x\in S$.
\end{cor}
\textbf{\emph{Proof.}} With Theorem (\ref{3.3}), if the set
$g(N_{g,\psi})$ is an unbounded set, then $f=T$, which implies
$f(x)=g(x)f(1)$ for all $x\in S$ and since
$$(g(x+y)-g(x)g(y))T(z)=0$$ for any $x,y,z\in S$ and $f=T$ is a
nonzero function for any nonzero function $f$, so $f\in
Z_{\Im(S,\mathbb{C})}$.\\

In \cite{Baker1}, Baker presented an example to show that
$$\|f(x+y)-f(x)f(y)\|\leq\varepsilon\ \ \ \ \ \ for\ x, y\in S$$
implies that $f$ is either bounded or exponential fails if the
algebra does not have the multiplicative norm property. Here, we
extend this result another way and conditions in unitary complex
Banach algebra.
\begin{thm}\label{3.10}
Let $S$ be commutative semigroup, let $B$ be a complex Banach
algebra with unit $1_{B}$, and let $\Im(S,B)$ be the functional
equation $f(x\cdot y)=g(x)f(y)$, where $g:S\rightarrow B$ is a
given function. If $N_{(\widehat{g},\psi)}\neq{\O}$, then the
functional equation $\Im(S,B)$ is HUR-stable on $(S,B)$ with
controls $(\psi(x,y),\inf_{a\in
N_{\widehat{g},\psi}}\frac{\psi(a,y)}{|\widehat{g}(a)|-1})$ such
that the HUR-stable function such $T$ is unique and
\begin{equation}\label{33.2}
(g(x\cdot y)-g(x)g(y))T(z)=0
\end{equation}
for all $x,y, z\in S$.
\end{thm}
\textbf{\emph{Proof.}} Let $a\cap N_{\widehat{g},\psi}$ be a
arbitrary fixed. We have
\begin{eqnarray*}
  \|f(a\cdot y)-g(a)f(y)\| &=& \|f(a\cdot y)-\widehat{g}(a)(1_{B} f(y))\| \\
   &=& \|f(a\cdot y)-\widehat{g}(a)f(y)\|\leq\psi(a,y)
\end{eqnarray*}
thus, $\|f(y\cdot a)-\hat{g}(a)f(y)\|\leq\psi(a,y)$ for all $y\in
S$. So from Theorem (\ref{3.3}), there is a unique function
$T:S\rightarrow B$ such that
$$T(x\cdot y)=\widehat{g}(x)T(y)$$
$$[\widehat{g}(x\cdot y)-\widehat{g}(x)\widehat{g}(y)]T(z)=0$$
and satisfying
\begin{equation}
\|f(y)-T(y)\|\leq\inf_{a\in
N_{\widehat{g}}}[\frac{\psi(a,y)}{|\widehat{g}(a)|-1}]
\end{equation}
for all $x, y, z\in S$. We have
\begin{equation}
  \|f(x\cdot y \cdot a^{n})-g(x)f(y\cdot a^{n})\|\leq\psi(x,y\cdot a^{n})
\end{equation}
then on the dividing by $|\widehat{g}(a)|^{n}$ we see that
\begin{equation}
\| \frac{f(x\cdot y \cdot
a^{n})}{\widehat{g}(a)^{n}}-g(x)\frac{f(y\cdot
a^{n})}{\widehat{g}(a)^{n}}\| \leq \frac{\psi(x,y\cdot
a^{n})}{|\widehat{g}(a)|^{n}} \leq
\frac{\psi(x,y)}{|\widehat{g}(a)|^{n}}.
\end{equation}
Hence, $T(x\cdot y)=g(x)T(y)$ for all $x, y\in S$. Now let $x, y,
z\in S$ be arbitrary elements, then
$$T(x\cdot y\cdot z)=g(x\cdot y)T(z)$$
and
$$T(x\cdot y\cdot z)=g(x)T(y\cdot z)=g(x)g(y)T(z)$$
so, $$(g(x\cdot y)-g(x)g(y))T(z)=0$$ for all $x, y, z\in S$. The
proof is complete.\\

In the following, we generalize the well-known Baker's
superstability and stability result for exponential mappings with
values in the field of complex numbers to the case of an arbitrary
unitary complex Banach algebra.
\begin{cor}
Let $(S,\cdot)$ be  commutative semigroup, $B$ be be a complex
Banach algebra with unit $1_{B}$, and let $\Im(S,B)$ be the
functional equation $f(x\cdot y)=f(x)f(y)$. If $f:S\rightarrow B$
a function such that $f(N_{\widehat{f},\psi})$ is an unbounded
set, then $f\in Z_{\Im(S,B)}$.

\end{cor}
\textbf{Proof.} In Theorem (\ref{3.10}), if we put $g(x):=f(x)$,
then we will had
\begin{center}
$\|f(x)-T(x)\|\leq\inf_{a\in
N_{\widehat{f},\psi}}[\frac{\psi(a,x)}{|f(a)|-1}]$ and
$(f(x+y)-f(x)f(y))T(z)=0$
\end{center}
for all $x, y, z\in S$. Now since
$\widehat{f}(N_{\widehat{f},\psi})$ is unbounded set, then we
have $f=T$, which says that $f\in Z_{\Im_{f}}$ and the proof is
complete.
\begin{cor}\label{homo4}
Let $(S,\cdot)$ be commutative semigroup with identity $1$, $B$
be a complex Banach algebra with unit $1_{B}$, and let $\Im(S,B)$
be the functional equation $f(1)f(x \cdot y)=g(x)f(y)$, where
$g:S\rightarrow \mathbb{C}$ is a given function. If
$\widehat{g}(N_{\widehat{g},\psi})$ is unbounded set, then the
functional equation $\Im(S,B)$ is hyperstable on $(S,B)$ with
control $(\psi(x,y))$ and $f(x)=g(x)f(1)$ for all $x\in S$.
\end{cor}
\textbf{\emph{Proof.}} With Theorem (\ref{3.10}).\\

In the following Theorem, we consider the superstability of the a
Pexiderized of exponential equation $$f(x+y)=g(x)h(y),$$ in which
$f, g$ and $h$ are three functions from a commutative semigroup
to an unitary an complex Banach algebra.
\begin{thm}\label{3.12}
Let $(S,\cdot)$ be commutative semigroup and $B$ be a complex
Banach algebra with unit $1_{B}$. Let $f,g,h:S\rightarrow B$ be
three functions and $g(x_{0})=1_{B}$ for a fixed  $x_{0}\in S$
and also
\begin{eqnarray}
  \|f(x\cdot y)-g(x)h(y)\| &\leq& \psi(x,y) \label{3.13}
\end{eqnarray}
for all $x,y\in S$. If $N_{\widehat{g},\psi}\neq{\O}$, then there
exists a exactly one function $T:S\rightarrow B$ such that
$$T(x\cdot y)=g(x)T(y),$$
$$(g(x\cdot y)-g(x)g(y))T(z)=0$$
and satisfies
$$\|f(y)-T(y)\|\leq\inf_{a\in N_{\widehat{g},\psi}}\frac{\widetilde{\psi}(a,y)}{|\widehat{g}(a)|-1},$$
$$\|h(y)-T(y)\|\leq\inf_{a\in N_{\widehat{g},\psi}}\frac{\widehat{\psi}(a,y)}{|\widehat{g}(a)|-1}$$
for all $x, y ,z\in S$, in which $\widetilde{\psi}(x,y)=\psi(x,y)+
\|g(x)\|\psi(x_{0},y)$ and
$\widehat{\psi}(x,y)=\psi(x,y)+\psi(x_{0},x\cdot y)$ for $x,y\in
S$.
\end{thm}
\textbf{Proof.} Applying (\ref{3.13}) we get for all $x,y\in S$
\begin{eqnarray*}
  \|f(x\cdot y)-g(x)f(y)\|&\leq& \|f(x\cdot y)-g(x)h(y)\|+\|g(x)f(y)-g(x)h(y)\| \\
   &\leq& \psi(x,y)+ \|g(x)\|\psi(x_{0},y)
\end{eqnarray*}
and
\begin{eqnarray*}
  \|h(x\cdot y)-g(x)h(y)\|&\leq& \|h(x\cdot y)-f(x\cdot y)\|+\|f(x\cdot y)-g(x)h(y)\| \\
   &\leq&  \psi(x,y)+|g(x_{0})|\psi(x_{0},x\cdot y)
\end{eqnarray*}
We set $\widetilde{\psi}(x,y)=\psi(x,y)+ \|g(x)\|\psi(x_{0},y)$
and $\widehat{\psi}(x,y)=\psi(x,y)+\psi(x_{0},x\cdot y)$ for
$x,y\in S$ and these are obvious that
$$\widetilde{\psi}(x,y\cdot a)\leq \widetilde{\psi}(x,y)$$ and
$$\widehat{\psi}(x,y\cdot a)\leq \widehat{\psi}(x,y)$$  for $x,y\in S$ and $a\in N_{\hat{g},\psi}$.
Therefore by Theorem (\ref{3.10}), then there exists a exactly one
function $H:S\rightarrow B$ such that
$$H(x+y)=g(x)H(y)$$
$$(g(x+y)-g(x)g(y))H(z)=0$$
and satisfies
$$\|f(y)-H(y)\|\leq\inf_{a\in N_{\widehat{g},\psi}}\frac{\widetilde{\psi}(a,y)}{|\widehat{g}(a)-1|}$$
for all $x, y ,z\in S$, where
$H(x)=\lim_{n\rightarrow\infty}\frac{f(x\cdot
a^{n})}{\widehat{g}(a)^{n}}$ for all $x\in S$ and any fixed $a\in
N_{\widehat{g},\psi}$. And also then there exists a exactly one
function $F:S\rightarrow B$ such that
$$F(x+y)=g(x)F(y)$$
$$(g(x+y)-g(x)g(y))F(z)=0$$
and satisfies
$$\|h(y)-F(y)\|\leq\inf_{a\in N_{\widehat{g},\psi}} \frac{\widehat{\psi}(a,y)}{|\widehat{g}(a)|-1}$$
for all $x, y ,z\in S$, where
$F(x)=\lim_{n\rightarrow\infty}\frac{h(x\cdot
a^{n})}{\widehat{g}(a)^{n}}$ for all $x\in S$ and any fixed $a\in
N_{\widehat{g},\psi}$. Furthermore, we have
$$\|\frac{f(x\cdot a^{n})}{\widehat{g}(a)^{n}}-\frac{h(x\cdot a^{n})}{\widehat{g}(a)^{n}}\|=|\widehat{g}(a)|^{-n}\|f(x\cdot a^{n})-h(x\cdot
a^{n})\|$$
    $$\leq \frac{|g(x_{0})|\psi(x_{0},x\cdot a^{n})}{\hat{g}(a)^{n}}\leq \frac{|g(x_{0})|\psi(x_{0},x)}{\hat{g}(a)^{n}}$$
for all $x\in S$ and any fixed $a\in N_{\widehat{g},\psi}$. Hence,
$H=F$ and so there exists an exactly one function $T:S\rightarrow
X$ such that
$$T(x\cdot y)=g(x)T(y),$$
$$(g(x\cdot y)-g(x)g(y))T(z)=0$$
and satisfies
$$\|f(y)-T(y)\|\leq\inf_{a\in N_{\widehat{g},\psi}} \frac{\widetilde{\psi}(a,y)}{|\widehat{g}(a)|-1},$$
$$\|h(y)-T(y)\|\leq\inf_{a\in N_{\widehat{g},\psi}} \frac{\widehat{\psi}(a,y)}{|\widehat{g}(a)|-1}$$
for all $x, y ,z\in S$. The proof is complete.
\begin{cor}\label{q0.1}
Let $(S,\cdot)$ be commutative semigroup and $B$ be  complex
field. Let $f,g,h:S\rightarrow \mathbb{C}$ be three functions and
$g(x_{0})=1$ for a fixed  $x_{0}\in S$ and also
\begin{eqnarray}
  \|f(x\cdot y)-g(x)h(y)\| &\leq& \psi(x,y) \label{3.13}
\end{eqnarray}
for all $x,y\in S$. If $g(N_{g,\psi})$ is an unbounded set, then
$f(x\cdot y)=g(x)h(y)$ for all $x,y\in S$, $f=h$, and
$f(x)=f(1)g(x)$ for all $x\in S$.
\end{cor}
\textbf{Proof.} Applying Theorem (\ref{3.12}), we get, there
exists a exactly one function $T:S\rightarrow \mathbb{C}$ such
that
$$T(x\cdot y)=g(x)T(y),$$
$$(g(x\cdot y)-g(x)g(y))T(z)=0$$
and satisfies
$$\|f(y)-T(y)\|\leq\inf_{a\in N_{\widehat{g},\psi}} \frac{\widetilde{\psi}(a,y)}{|\widehat{g}(a)|-1},$$
$$\|h(y)-T(y)\|\leq\inf_{a\in N_{\widehat{g},\psi}} \frac{\widehat{\psi}(a,y)}{|\widehat{g}(a)|-1}$$
for all $x, y ,z\in S$. Since $g=\widehat{g}$ and $g(N_{g,\psi})$
is an unbounded set, so $f=h=T$ and so that $f(x)=g(x)f(1)$ for
all $x\in S$. The proof is complete.\\

Therefore, with above Corollary, we give a version of Baker's
Theorem (\ref{b1}) for Pexiderized of exponential functional
equation.
\begin{cor}\label{homo5}
Let $(S,\cdot)$ be commutative semigroup and $B$ be complex
field. Let $f,g,h:S\rightarrow \mathbb{C}$ be three functions and
$g(x_{0})=1$ for a fixed  $x_{0}\in S$ and also
\begin{eqnarray}
  \|f(x\cdot y)-g(x)h(y)\| &\leq& \epsilon \label{3.13}
\end{eqnarray}
for all $x,y\in S$ for some $\epsilon>0$. If the function $g$ is
an unbounded function, then $f(x\cdot y)=g(x)h(y)$ for all
$x,y\in S$, $f=h$, and $f(x)=f(1)g(x)$ for all $x\in S$.
\end{cor}
\textbf{Proof.} If we set $\psi(x,y):=\epsilon$ for all $x,y \in
S$. Its obvius that $N_{g,\epsilon}=\{a\in S\ :\ |g(a)|>1\}$,
where implies that $g(N_{g,\epsilon})$ is an unbounded set and so
that Corollary (\ref{q0.1}) complete the proof.

\begin{center}{\textbf{Some Remarks About Homogeneous Functional Equation}}
\end{center}
The functional equation
\begin{equation}\label{hom1}
  f(yx)=y^{k}f(x)
\end{equation}
(where $k$ is a fixed real constant) is called the "homogeneous
functional equation" of degree $k$. In the case when $k=1$ in the
equation (\ref{hom1}), the equation is simply called the
"homogeneous functional equation". In general, "$\phi$-homogeneous
functional equation" is the following equation:
\begin{equation}\label{hom2}
  f(\alpha x)=\phi(\alpha)f(x)
\end{equation}
for suitable class function and spaces. In \cite{Czerwik1}, S.
Czerwik considered the Pexiderized homogeneous functional
equation and he obtained the following result:
\begin{thm}
Let $V$ be a real linear space and $E$ a real Banach space. Let
$f,g:V\rightarrow E$ and $\phi: \mathbb{R}\rightarrow \mathbb{R}$
and $h: \mathbb{R}\times E\rightarrow \mathbb{R}^{+}$ be given
functions. Assume that
$$\|f(\alpha x)-\phi(\alpha)g(x)\|\leq h(\alpha,x)$$
for all $(\alpha,x)\in \mathbb{R}\times V$ and $\phi(1)=1$.
Suppose that there exists $\beta\in \mathbb{R}$ such that
$\phi(\beta)\neq 0$ and the series
$$\sum_{n=1}^{\infty}|\phi(\beta)|^{-n}H(\beta,\beta^{n}x)$$
converges pointwise for all $x\in V$, and
$$\liminf_{n\rightarrow\infty}|\phi(\beta)|^{-n}H(\alpha,\beta^{n}x)=0$$
for all $(\alpha,x)\in \mathbb{R}\times V$, where $H(\alpha,x):=
h(\alpha,x)+|\phi(\alpha)|h(1,x)$. Then there exists exactly one
$\phi$-homogeneous function $A:V\rightarrow E$:
$$A(\alpha x)=\phi(\alpha)A(x)$$
for all $(\alpha,x)\in \mathbb{R}\times V$ such that
$$\|f(x)-A(x)\|\leq \sum_{n=1}^{\infty}|\phi(\beta)|^{-n}H(\beta,\beta^{n-1}x),$$
$$\|g(x)-A(x)\|\leq \sum_{n=1}^{\infty}|\phi(\beta)|^{-n}G(\beta,\beta^{n-1}x),$$
for all $x\in V$, where $G(\alpha,x):= h(\alpha,x)+h(1,\alpha x)$.
\end{thm}
For more infomation about this subject one can refer to
(\cite{Jung1}, Chap. 5). Since in this section, we study the
exponential funcyional equation and its Pexiderized on seimigroup
domain, so the results (\ref{3.3}), (\ref{homo3}), (\ref{3.10}),
(\ref{homo4}), (\ref{3.12}), (\ref{q0.1}), and (\ref{homo5}) is
applied for study of superstability and stability of the
homogeneous functional equation and its Pexiderized. Hence, we
don't present the similarly results for homogeneous functional
equation.
\section{\bf Superstability of the Cauchy Additive Functional Equation on Semigroups}
The results of this section was proved by Mohsen Alimohammady and
Ali Sadeghi, which was published see \cite{sadeghi2}. Throughout
of this section, assume that $(S,\cdot)$ is an arbitrary
commutative semigroup, $\mathbb{C}$ is the field of all complex
numbers, $\mathbb{R}$ is real field, $\widehat{\mathbb{C}}=\{z\in
\mathbb{C}\ |\ -\pi<\arg(z)\leq\pi\}$, $\psi:S\rightarrow
\mathbb{R}^{+}$ and $\phi:S^{2}\rightarrow \mathbb{R}^{+}$ are
some functions. Also for the function $f:S\rightarrow
\widehat{\mathbb{C}}$, then $f(S)^{+}$ is a subset of
$\widehat{\mathbb{C}}$, where $f(S)^{+}=\{p\in S\ |\
Re(f(p))>0\}$. In this section, we call $f:S\rightarrow
\widehat{\mathbb{C}}$ is a Cauchy function, if
$$f(x\cdot y)-f(x)-f(y)=0$$ for all $x,y\in S$.
\begin{thm}\label{4.1}
Let $f:S\rightarrow \widehat{\mathbb{C}}$ is a function and
\begin{eqnarray}
  |f(x\cdot y)-f(x)-f(y)|&\leq&\phi(x,y) \label{4.2};\\
  |f(x)|&\leq& \psi(x) \label{4.3}
\end{eqnarray}
for all $x, y\in S$. Assume that there exists $p\in f(S)^{+}$ such
that
\begin{eqnarray}
  \sum_{m=0}^{\infty}\phi(p,p^{m+1})<\infty \label{4.4}; \\
  \psi(x\cdot p)&\leq&\psi(x)
\end{eqnarray}
for all $x\in S$. Then $f$ is a Cauchy function.
\end{thm}
Note that the above Theorem is partial affirmative answer to
problem $18$, in the thirty-first ISFE. Moreover, we present a
superstability result for the logarithmic functional equation.

\emph{\textbf{Proof of Theorem.}} Let
$E:\mathbb{C}\rightarrow\mathbb{C}$ be exponential function, where
$E(a)=\exp(a)$ for each $a\in \mathbb{C}$. Now from (\ref{4.3}),
we have
\begin{eqnarray*}
  |E(f(x\cdot y))-E(f(x)+f(y))| &\leq& |E(f(x\cdot y))|+|E(f(x)+f(y))| \\
    &\leq& E(|f(x\cdot y)|)+E(|f(x)+f(y)|) \\
    &\leq& E(\psi(x\cdot y))+E(\psi(x)+\psi(y))
\end{eqnarray*}
for all $x,y\in S$. So, the function
$\widehat{E}:S\rightarrow\mathbb{C}$ with
$\widehat{E}=(Eof)(x)=\exp(f(x))$ for all $x\in S$, satisfies the
following inequality
\begin{eqnarray}\label{4.5}
 |\widehat{E}(x\cdot y)-\widehat{E}(x)\widehat{E}(y)| &\leq& \varphi(x,y)
\end{eqnarray}
for all $x,y\in S$, in which $\varphi(x,y):= E(\psi(x\cdot
y))+E(\psi(x)+\psi(y))$ for all $x, y\in S$. From assumption
$\psi(x\cdot p)\leq\psi(x)$ for all $x\in S$. So, its easy to show
that
\begin{eqnarray}\label{4.6}
 \varphi(x,y\cdot p)&\leq& \varphi(x,y)
\end{eqnarray}
for all $x,y\in S$. We are going to show that $\widehat{E}(x\cdot
y)-\widehat{E}(x)\widehat{E}(y)=0$ for all $x, y\in S$. Now let us
consider the set $\mathcal{A}:=\{h:S\rightarrow \mathbb{C}\}$ and
introduce the generalized metric on $\mathcal{A}$:
$$d(u,h)=\sup_{x\in S}\frac{|u(x)-h(x)|}{\varphi(p,x)}.$$
It is easy to show that $(\mathcal{A}, d)$ is generalized
complete metric space. Now we define the function
$J:\mathcal{A}\rightarrow \mathcal{A}$ with
$$J(h(x))=\frac{1}{\widehat{E}(p)}h(p.x)$$ for all $h\in \mathcal{A}$ and $x\in S$. Since $\varphi(x,y\cdot p)\leq\varphi(x,y) $ for all
$x\in S$, so
  \begin{eqnarray*}
    d(J(u),J(h)) &=& \sup_{x\in S}\frac{|u(p.x)-h(p.x)|}{|\widehat{E}(p)|\varphi(p,x)} \\
      &\leq& \sup_{x\in S}\frac{|u(p.x)-h(p.x)|}{|\widehat{E}(p)|\varphi(p,p\cdot x)} \\
      &\leq& \frac{1}{|\widehat{E}(p)|}d(u,h)
  \end{eqnarray*}
for all $u, h\in \mathcal{A}$, that is $J$ is a strictly
contractive selfmapping of $\mathcal{A}$, with the Lipschitz
constant $L=\frac{1}{|\widehat{E}(p)|}$ (note that
$|\widehat{E}(p)|=|exp(f(p))|=exp(Re(f(p)))>1$). From (\ref{4.5}),
we get
$$|\frac{\widehat{E}(p\cdot x)}{\widehat{E}(p)}-\widehat{E}(x)|\leq\frac{\varphi(p,x)}{|\widehat{E}(p)|}$$
for all $x\in S$, which says that
$d(J(\widehat{E}),\widehat{E})\leq L<\infty$. By Theorem
(\ref{3.1}), there exists a mapping $T:G\rightarrow \mathbb{C}$
such that
  \begin{enumerate}
  \item $T$ is a unique fixed point of $J$, i.e.,
  \begin{eqnarray}
    T(p\cdot x) &=& \widehat{E}(p)T(x)
  \end{eqnarray}
  for all $x\in S$.
  \item $d(J^{n}(\widehat{E}),T)\rightarrow 0$ as $n\rightarrow \infty$. This implies that
  \begin{eqnarray}\label{4.7}
    T(x)&=&lim_{n\rightarrow\infty}\frac{\widehat{E}(p^{n}\cdot x)}{\widehat{E}^{n}(p)}
  \end{eqnarray}
  for all $x\in S$.
  \item $d(\widehat{E},T)\leq\frac{1}{1-L}d(J(\widehat{E}),\widehat{E})$, which implies,
  $$d(\widehat{E},T)\leq\frac{1}{|\widehat{E}(p)|-1}.$$
 \end{enumerate}
Let $x,y\in S$  be two arbitrary fixed elements, from (\ref{4.5})
and (\ref{4.6})
$$|\widehat{E}(p^{n}\cdot(x\cdot y))-\widehat{E}(x)\widehat{E}(p^{n}\cdot y)|\leq\varphi(x,p^{n}\cdot y)$$ and dividing by $|\widehat{E}(p)|^{n}$,
$$|\frac{\widehat{E}(p^{n}\cdot(x\cdot y))}{\widehat{E}(p)^{n}}-\widehat{E}(x)\frac{\widehat{E}(p^{n}\cdot y)}{\widehat{E}(p)^{n}}|\leq
\frac{\varphi(x,p^{n}\cdot
y)}{|\widehat{E}(p)|^{n}}\leq\frac{\varphi(x,y)}{|\widehat{E}(p)|^{n}}$$
and letting $n$ to infinity, we get $T(x\cdot
y)=\widehat{E}(x)T(y)$, which says that $T$ is a
$\widehat{E}$-homogeneous function.

Let $x,y,z \in S$ be arbitrary elements, then
$$T((x\cdot y)\cdot z)=\widehat{E}(x\cdot y)T(z)$$ and $$T((x\cdot y)\cdot z)=\widehat{E}(x)T(y\cdot z)=\widehat{E}(x)\widehat{E}(y)T(z)$$
or
\begin{eqnarray}\label{4.8}
  (\widehat{E}(x\cdot y)-\widehat{E}(x)\widehat{E}(y))T(z)&=&0.
\end{eqnarray}
Now we show that $T(p)\neq0$. From (\ref{4.2}), its easy to show
that following inequality
$$|Re(f(p^{n+1})-(n+1)f(p))|\leq|f(p^{n+1})-(n+1)f(p)|<\sum_{i=0}^{n-1}\phi(p,p^{i+1})$$
for all $n\in \mathbb{N}$. Also, from relation (\ref{4.7}), we
obtain
\begin{eqnarray}
  |T(p)| &=& \lim_{n\rightarrow\infty}|\frac{\widehat{E}(p^{n+1})}{\widehat{E}^{n}}| \\
   &=& \lim_{n\rightarrow\infty} |exp(f(p^{n+1})-(n+1)f(p))| \\
   &=& \lim_{n\rightarrow\infty} exp(Re(f(p^{n})-(n)f(p))).
\end{eqnarray}
Since $|T(p)|<\infty$, So $\lim_{n\rightarrow\infty}
exp(Re(f(p^{n})-(n)f(p)))$ there exist. Now from (\ref{4.6}) and
(\ref{4.4}), we obtain
$$\lim_{n\rightarrow \infty}|Re(f(p^{n+1})-(n+1)f(p))|<\infty,$$
where it implies that $|T(p)|\neq0$ . From (\ref{4.8}), we get
$$\widehat{E}(x\cdot y)=\widehat{E}(x)\widehat{E}(y)$$
or
$$\exp(f(x+y))=\exp(f(x)+f(y))$$
for all $x,y\in S$. Since exponential function is one-to-one on
$\widehat{\mathbb{C}}$, so $f$ is a Cauchy function. The proof is
complete.
\begin{cor}\label{4.10}
Let $S$ be a multiplicative semigroup of $\mathbb{R}$ and
$L:S\rightarrow \mathbb{R}$ be a function such that
\begin{eqnarray}
  |L(xy)-L(x)-L(y)|&\leq&\phi(x,y) ;\\
  |L(x)|&\leq& \psi(x)
\end{eqnarray}
for all $x, y\in S$. Assume that there exists $p\in S$ such that
$L(p)\neq0$ and
\begin{eqnarray}
  \sum_{m=0}^{\infty}\phi(p,p^{m+1})<\infty ; \\
  \psi(xp)&\leq&\psi(x)
\end{eqnarray}
for all $x\in S$. Then $L$ is a Cauchy function or a logarithmic
function i.e.; $L(xy)=L(x)+L(y)$ for all $x,y\in S$.
\end{cor}
\emph{\textbf{Proof.}} If $L(p)<0$, set $f:=-L$ and if $L(p)>0$,
set $f:=L$. Now applying Theorem (\ref{4.1}), we get the result.
\section{\bf Asymptotic Behavior of Cauchy Additive Functional Equations}

The results of this section was proved by Mohsen Alimohammady and
Ali Sadeghi, which was published see \cite{sadeghi3}. Throughout
this section, assume that $(S,+)$ is a commutative semigroup,
$E_{1}$, $E_{2}$ are two complex normed spaces, $\mathbb{R}$ is
real field, $\mathbb{N}$ is all positive integers, and
$\psi:S^{2}\rightarrow [0,\infty)$ is a function.

The following Theorem is a affirmative answer to problem $18$, in
the thirty-first ISFE.
\begin{thm}\label{5.1}
Let $\Im(S,E_{1})$ be the functional equation $f(x+y)=f(x)+f(y)$.
Assume that
\begin{itemize}
  \item $\lim_{n\rightarrow \infty}\frac{1}{n}\sum_{i=0}^{n-1}\psi(x+ix_{0},x_{0})=0$;
  \item $\lim_{n\rightarrow \infty}\frac{1}{n}\psi(x+nx_{0},y+ny_{0})=0$
\end{itemize}
for any fixed $x_{0}, y_{0}, x, y\in S$. Then the functional
equation $\Im(S,E_{1})$ is hyperstable on $(S,E_{1})$ with control
$(\psi)$.
\end{thm}
\textbf{\emph{Proof.}} Let $x_{0}$ be any fixed element of $S$
and we have
$$\|f(x+x_{0})-f(x_{0})-f(x)\|\leq\psi(x,x_{0})$$
for all $x\in S$. From last inequality, its easy to show that the
following inequality
$$\|f(x+nx_{0})-nf(x_{0})-f(x)\|\leq\sum_{i=0}^{n-1}\psi(x+ix_{0},x_{0})$$
for each fixed $x\in S$ and $n\in \mathbb{N}$. Now bye assumption
$\lim_{n\rightarrow
\infty}\frac{1}{n}\sum_{i=0}^{n-1}\psi(x+ix_{0},x_{0})=0$, so
$$f(x_{0})=\lim_{n\rightarrow \infty}\frac{f(x+nx_{0})}{n}$$
for any fixed $x\in S$. Let $x_{0}, y_{0}$ be any two fixed
elements of $S$, then we have
$$\|f(x+y+n(x_{0}+y_{0}))-f(x+nx_{0})-f(y+ny_{0})\|\leq\psi(x+nx_{0},y+ny_{0})$$
for any fixed $x, y\in S$. Now since $\lim_{n\rightarrow
\infty}\frac{1}{n}\psi(x+nx_{0},y+ny_{0})=0$, thus
$$f(x_{0}+y_{0})=f(x_{0})+f(y_{0}),$$
which says that $f\in Z_{\Im(S,B)}$ and the proof is complete.\\

As a consequence of the following result, its show that the
Theorem of Rassias (\ref{Rassias}) for the case $p<0$ is
hypestable.
\begin{cor}\label{2r}
Let $\Im(E_{1},E_{2})$ be the functional equation
$f(x+y)=f(x)+f(y)$. Then for every real's $\epsilon>0$, $p<0$ and
$q\leq 1$, the functional equation $\Im(S,B)$ is hyperstable on
$(S,B)$ with control $(\epsilon(\|x\|^{p}+\|y\|^{q}))$.
\end{cor}
\textbf{\emph{Proof.}} Set $\psi(x,y):=
(\epsilon(\|x\|^{p}+\|y\|^{q})$ for all $x, y\in S$. Since the
sequence $\sum_{i=0}^{n-1}\psi(x+ix_{0},x_{0})$ for any fixed $x,
y, x_{0}, y_{0}$ is increasing sequances and also $p<0$ and
$q\leq 1$, so that obvious that the followings relations:
\begin{itemize}
  \item $\lim_{n\rightarrow \infty}\frac{1}{n}\sum_{i=0}^{n-1}\psi(x+ix_{0},x_{0})=0$;
  \item $\lim_{n\rightarrow \infty}\frac{1}{n}\psi(x+nx_{0},y+ny_{0})=0$
\end{itemize}
for any fixed $x_{0}, y_{0}, x, y\in S$. Now Theorem (\ref{5.1})
implies that the result and the proof is complete.\\

In the following, by using Theorem (\ref{5.1}), we give a simple
proof of Skof theorem (\ref{11}) and also we show that Skof
theorem is true when $E_{2}$ be a complex normed space.
\begin{thm}
For a function $f:E_{1}\rightarrow E_{2}$ the following two
conditions are equivalent:
\begin{enumerate}
  \item $\|f(x+y)-f(x)-f(y)\|\rightarrow 0$ as $\|x\|+\|y\|\rightarrow \infty$;
  \item $f(x+y)-f(x)-f(y)=0$
\end{enumerate}
for all $x, y\in E_{1}$.
\end{thm}
\textbf{\emph{Proof.}} Set $\psi(x,y):=\|f(x+y)-f(x)-f(y)\|$ for
all $x,y\in E_{1}$. Now let $x_{0}, y_{0}\in E_{1}$ be two
arbirary fixed elements. Since
$\|x+nx_{0}\|+\|y+ny_{0}\|\rightarrow \infty$ for each fixed $x,
y\in E_{1}$, so
$$\lim_{n\rightarrow \infty}\psi(x+nx_{0},y+ny_{0})=0,$$
for each fixed $x, y\in E_{1}$, hence its easy to show that the
following relations
\begin{itemize}
  \item $\lim_{n\rightarrow \infty}\frac{1}{n}\sum_{i=0}^{n-1}\psi(x+ix_{0},x_{0})=0$;
  \item $\lim_{n\rightarrow \infty}\frac{1}{n}\psi(x+nx_{0},y+ny_{0})=0$
\end{itemize}
for each fixed $x, y\in E_{1}$. Now with Theorem (\ref{5.1})
implies that $f$ is an additive mapping (i.e.,
$f(x+y)=f(x)+f(y)$). The proof is complete.\\

Let $\mathfrak{S}$ be set all function $\rho:
E_{1}^{2}\rightarrow [0,\infty)$ such that
\begin{enumerate}
  \item $\rho(x+nx_{0},y+ny_{0})\rightarrow \infty$ as $n\rightarrow \infty$
\end{enumerate}
for any fixed $x_{0}, y_{0}, x, y\in E_{1}$. Not that the
functions $\rho_{1},\rho_{2},\rho_{3}\in \mathfrak{S}$, in which
$\rho_{1}(x,y):=\|x\|+\|y\|$, $\rho_{2}(x,y):=\|x+y\|$ and
$\rho_{3}(x,y):=\max\{\|x\|,\|y\|\}$ for all $x, y\in E_{1}$. We
now apply Theorem (\ref{5.1}) to a generalization of Skof theorem.
\begin{cor}
For a function $f:E_{1}\rightarrow E_{2}$ the following two
conditions are equivalent:
\begin{enumerate}
  \item $\|f(x+y)-f(x)-f(y)\|\rightarrow 0$ as $\rho(x,y)\rightarrow \infty$;
  \item $f(x+y)-f(x)-f(y)=0$
\end{enumerate}
for all $x, y\in E_{1}$, in which $\rho\in \mathfrak{S}$.
\end{cor}
\textbf{\emph{Proof.}} Set $\psi(x,y):=\|f(x+y)-f(x)-f(y)\|$ for
all $x,y\in E_{1}$. Now let $x_{0}, y_{0}\in E_{1}$ be two
arbirary fixed elements. Since $\rho\in \mathfrak{S}$, so
  $$\rho(x+nx_{0},y+ny_{0})\rightarrow \infty$$
as $n\rightarrow \infty$ for any fixed $x_{0}, y_{0}, x, y\in
E_{1}$. Thus
$$\lim_{n\rightarrow \infty}\psi(x+nx_{0},y+ny_{0})=0$$
for each fixed $x_{0}, y_{0}, x, y\in E_{1}$. Hence, its easy to
show that the following relations:
\begin{itemize}
  \item $\lim_{n\rightarrow \infty}\frac{1}{n}\sum_{i=0}^{n-1}\psi(x+ix_{0},x_{0})=0$;
  \item $\lim_{n\rightarrow \infty}\frac{1}{n}\psi(x+nx_{0},y+ny_{0})=0$
\end{itemize}
for any fixed $x_{0}, y_{0}, x, y\in E_{1}$. Now with Theorem
(\ref{5.1}) implies that $f$ is an additive mapping (i.e.,
$f(x+y)=f(x)+f(y)$). The proof is complete.
\subsection{\bf Asymptotic Behavior of Pexiderized Cauchy Additive Functional Equation}
\begin{thm}\label{5.5}
Let $S$ be with identity $e$ and $f,g,h:S\rightarrow E_{1}$ be
three functions such that $g(e)=h(e)=0$ and
\begin{eqnarray}\label{5.6}
  \|f(x+y)-g(x)-h(y)\|\leq\psi(x,y)
\end{eqnarray}
for all $x,y\in S$. Assume that
\begin{itemize}
  \item $\lim_{n\rightarrow \infty}\frac{1}{n}\sum_{i=0}^{n-1}\psi(x+ix_{0},x_{0})=0$;
  \item $\lim_{n\rightarrow \infty}\frac{1}{n}\psi(x+nx_{0},y+ny_{0})=0$
\end{itemize}
for any fixed $x_{0}, y_{0}, x, y \in S$. Then $f$, $g$ and $h$
are additive function and $f(x+y)-g(x)-h(y)=0$ for all $x, y\in
S$.
\end{thm}
\textbf{\emph{Proof.}} Set
$\widetilde{\psi}(x,y):=\psi(x,y)+\psi(x,e)+\psi(e,y)$ and
$\widehat{\psi}(x,y):=\psi(x+y,e)+\psi(x,e)+\psi(e,y)$ for all
$x, y \in S$.  From inequality (\ref{5.6}) and assumptions, we
obtain the following inequalities
\begin{eqnarray*}
  \|f(x+y)-f(x)-f(y)\| &\leq& \psi(x,y)+\|f(x)-g(x)\|+\|f(y)-h(y)\| \\
    &\leq& \psi(x,y)+\psi(x,e)+\psi(e,y)=\widetilde{\psi}(x,y)
\end{eqnarray*}
and
\begin{eqnarray*}
  \|g(x+y)-g(x)-g(y)\| &\leq& \psi(x+y,e)+\|f(x+y)-g(x)-g(y)\|\\
    &\leq& \widetilde{\psi}(x,y)+\widehat{\psi}(x,y)
\end{eqnarray*}
and also
\begin{eqnarray*}
  \|h(x+y)-h(x)-h(y)\| &\leq& \psi(x+y,e)+\|f(x+y)-h(x)-h(y)\|\\
    &\leq& \widetilde{\psi}(x,y)+\widehat{\psi}(x,y)
\end{eqnarray*}
for all $x, y\in S$. With assumptions its easy to show that
\begin{itemize}
  \item $\lim_{n\rightarrow \infty}\frac{1}{n}\sum_{i=0}^{n-1}\phi(x+ix_{0},x_{0})=0$;
  \item $\lim_{n\rightarrow \infty}\frac{1}{n}\phi(x+nx_{0},y+ny_{0})=0$
\end{itemize}
for any fixed $x_{0}, y_{0}, x, y \in S$, in which the function
$\phi$ is $\widetilde{\psi}$ or $\widetilde{\psi}+\widehat{\psi}$
. Now with Theorem $(\ref{5.1})$ $f, g$ and $h$ is additive
mapping and also
\begin{itemize}
  \item $f(x_{0})=\lim_{n\rightarrow \infty}\frac{f(x+nx_{0})}{n}$
  \item $g(x_{0})=\lim_{n\rightarrow \infty}\frac{g(x+nx_{0})}{n}$
  \item $h(x_{0})=\lim_{n\rightarrow \infty}\frac{h(x+nx_{0})}{n}$
\end{itemize}
for each fixed $x_{0}, x\in S$. Let $x_{0}, y_{0}$ be any two
fixed element of $S$, then from (\ref{5.6}), we obtain
$$\|f(x+y+n(x_{0}+y_{0}))-g(x+nx_{0})-h(y+ny_{0})\|\leq\psi(x+nx_{0},y+ny_{0})$$
for any fixed $x, y\in S$. Now since $\lim_{n\rightarrow
\infty}\frac{1}{n}\psi(x+nx_{0},y+y_{0})=0$, thus
$$f(x_{0}+y_{0})=g(x_{0})+h(y_{0}),$$
which says that  $f(x+y)-g(x)-h(y)=0$ for all $x, y\in S$. The
proof is complete.\\

In the following, by using Theorem (\ref{5.5}), we give a
generalization of Skof theorem for Pexiderized additive mapping.
\begin{thm}\label{5.7}
Assume that $f, g, h:E_{1}\rightarrow E_{2}$ are three functions
such that $g(0)=h(0)=0$, then the following two conditions are
equivalent:
\begin{enumerate}
  \item $\|f(x+y)-g(x)-h(y)\|\rightarrow 0$ as $\rho(x,y)\rightarrow \infty$;
  \item $f(x+y)-g(x)-h(y)=0$
\end{enumerate}
for all $x, y\in E_{1}$, in which $\rho\in \mathfrak{S}$.
\end{thm}
\textbf{\emph{Proof.}} Set $\psi(x,y):=\|f(x+y)-g(x)-h(y)\|$ for
all $x,y\in E_{1}$. Now let $x_{0}, y_{0}\in E_{1}$ be two
arbirary fixed elements. Since $\rho\in \mathfrak{S}$, so
$$\rho(x+nx_{0},y+ny_{0})\rightarrow \infty$$
as $n\rightarrow \infty$ for any fixed $x_{0}, y_{0}, x, y\in
E_{1}$. Thus
$$\lim_{n\rightarrow \infty}\psi(x+nx_{0},y+ny_{0})=0$$
for each fixed $x_{0}, y_{0}, x, y\in E_{1}$. Hence, its easy to
show that the following relations
\begin{itemize}
  \item $\lim_{n\rightarrow \infty}\frac{1}{n}\sum_{i=0}^{n-1}\psi(x+ix_{0},x_{0})=0$;
  \item $\lim_{n\rightarrow \infty}\frac{1}{n}\psi(x+nx_{0},y+ny_{0})=0$
\end{itemize}
for any fixed $x_{0}, y_{0}, x, y\in E_{1}$. Now with Theorem
(\ref{5.5}) implies that $f(x+y)-g(x)-h(y)=0$ for all $x, y\in
S$. The proof is complete.\\

In the following, with using Theorem (\ref{5.7}), we give a simple
proof of S.-M. Joung's theorem (see \cite{Joung}) and also we show
that S.-M. Joung's theorem is true when $E_{2}$ be a complex
normed space.
\begin{thm}\label{5.8}
Assume that $J:E_{1}\rightarrow E_{2}$ is a function such that
$J(0)=0$, then the following two conditions are equivalent:
\begin{enumerate}
  \item $\|2J(\frac{x+y}{2})-J(x)-J(y)\|\rightarrow 0$ as $\|x\|+\|y\|\rightarrow \infty$;
  \item $2J(\frac{x+y}{2})-J(x)-J(y)=0$
\end{enumerate}
for all $x, y\in E_{1}$.
\end{thm}
\textbf{\emph{Proof.}} Sets $f(x):=2J(\frac{x}{2})$, $g(x):=J(x)$,
and $\rho(x,y):=\|x\|+\|y\|$ for all $x, y\in E_{1}$. Now apply
Theorem (\ref{5.7}).
\vskip 0.4 true cm

\begin{center}{\textbf{Some Remarks}}
\end{center}
In 2013, J. Brzd\c{e}k \cite{Brzdek} proved that the Theorem of
Rassias (\ref{Rassias}) for the case $p<0$ is hypestable. Where
we prove this result earlier than J. Brzd\c{e}k and also we prove
its generalization (\ref{2r}), which was published in 2012 see
\cite{sadeghi3}. And also in 2013, M. Piszczek \cite{Piszczek}
consider a general Cauchy additive functional equation as follows
\begin{equation}\label{1r}
  g(ax+by)=Ag(x)+Bg(y)
\end{equation}
for class functions $g:X\rightarrow Y$, where $X$ is a normed
space over field $\mathbb{F}$, $Y$ is Banach space over
$\mathbb{F}$, and $\mathbb{F}$ is the fields of real or complex
numbers. He prove that the Theorem of Rassias (\ref{Rassias}) for
the case $p<0$ for functional equation (\ref{1r}) is hypestable.
Where this result is as consequence of Theorem (\ref{5.5}), and
also we consider more general calss functions, which was
published in 2012 see \cite{sadeghi3}.
\section{\bf Stability and Common Stability for the Systems of Linear Equations}

The results of this section was proved by Mohsen Alimohammady and
Ali Sadeghi, which was published see \cite{sadeghi4}. In this
section, First we consider the Hyers-Ulam-Rassias stability via a
fixed point approach for the linear functional equation
(\ref{e1}) and then applying these result we will investigate
Pexiderized linear functional equation (\ref{e4}).

Assume that $S$ is a nonempty set, $F=\mathbb{Q}$, $\mathbb{R}$
or $\mathbb{C}$, $B$ is a Banach spaces over $F$,
$\psi:S\rightarrow \mathbb{R}^{+}$, $f,g:S\rightarrow B$,
$p:S\rightarrow K\backslash\{0\}$, $q:S\rightarrow B$ are
functions, and $\sigma:S\rightarrow S$ is a arbitrary map.
\begin{thm}\label{6.1}
Let $\Im(S,B)$ be functional equation $f(\rho(x))=p(x)f(x)+q(x)$.
If there exists a real $0< L < 1$ such that
\begin{eqnarray}
  \psi(\rho(x))&\leq& L|p(\rho(x))|\psi(x)
\end{eqnarray}
for all $x\in S$. Then the functional equation $\Im(S,B)$ is
HUR-stable on $(S,B)$ with controls $(\psi(x),
\frac{\psi(x)}{(1-L)|p(x)|})$ and HUR-stable function such $T$ is
unique.
\end{thm}
\emph{\textbf{Proof.}} Let us consider the set
$\mathcal{A}:=\{h:S\rightarrow B\}$ and introduce the generalized
metric on $\mathcal{A}$:
$$d(u,h)=\sup_{\{x\in S\ ;\ \psi(x)\neq0\}}\frac{|p(x)|\|g(x)-h(x)\|}{\psi(x)}.$$
It is easy to show that $(\mathcal{A}, d)$ is generalized complete
metric space. Now we define the function $J
:\mathcal{A}\rightarrow \mathcal{A}$ with
$$J(h(x))=\frac{1}{p(x)}h(\rho(x))-\frac{q(x)}{p(x)}$$
for all $h\in \mathcal{A}$ and $x\in S$. Since $\psi(\rho(x))\leq
L|p(\rho(x))|\psi(x)$ for all $x\in X$, so
  \begin{eqnarray*}
    d(J(u),J(h)) &=& \sup_{\{x\in X\ ;\ \psi(x)\neq0\}}\frac{|p(x)|\|u(\rho(x))-h(\rho(x))\|}{|p(x)|\psi(x)} \\
      &\leq& \sup_{\{x\in X\ ;\ \psi(\rho(x))\neq0\}}L\frac{|p(\rho(x))|\|u(\rho(x))-h(\rho(x))\|}{\psi(\rho(x))}\leq Ld(u,h)
  \end{eqnarray*}
for all $u, h\in \mathcal{A}$, that is $J$ is a strictly
contractive selfmapping of $\mathcal{A}$, with the Lipschitz
constant $L$ (note that $0<L<1$). We have
$$\|f(\rho(x))-p(x)f(x)-q(x)\|\leq \psi(x)$$
for all $x\in S$, we get
$$\|\frac{f(\rho(x))}{p(x)}-\frac{q(x)}{p(x)}-f(x)\|\leq\frac{\psi(x)}{|p(x)|}$$
for all $x\in S$, which says that $d(J(f),f)\leq 1<\infty$. So,
with Theorem (\ref{3.1}), there exists a mapping $T:X\rightarrow
B$ such that
  \begin{enumerate}
  \item $T$ is a fixed point of $J$, i.e.,
  \begin{eqnarray}\label{6.3}
    T(\rho(x)) &=& p(x)T(x)+q(x)
  \end{eqnarray}
  for all $x\in S$. The mapping $T$ is a unique fixed point of $J$ in the set $\tilde{\mathcal{A}}=\{h\in \mathcal{A}\ :\ d(f,h)<\infty\}$.
  This implies that $T\in Z_{\Im(S,B)}$ and is unique HUR-stable
  function. Also there exists $C\in (0,\infty)$ satisfying
  $$\|f(x)-T(x)\|\leq C\frac{\psi(x)}{|p(x)|}$$
  for all $x\in X$.
  \item $d(J^{n}(f),T)\rightarrow 0$ as $n\rightarrow \infty$. This implies that
  $$T(x)=\lim_{n\rightarrow\infty}\frac{f(\rho^{n}(x))}{\prod_{i=0}^{n-1}p(\rho^{i}(x))}-\sum_{k=0}^{n-1}\frac{q(\rho^{i}(x))}{\prod_{i=0}^{k}p(\rho^{i}(x))}$$
  for all $x\in X$.
  \item $d(f,T)\leq\frac{1}{1-L}d(J(f),f)$, which implies,
  $$d(f,T)\leq\frac{1}{1-L}$$
  or
  $$\|f(x)-T(x)\|\leq \frac{\psi(x)}{(1-L)|p(x)|}$$
for all $x\in X$.
\end{enumerate}
Therefore, the functional equation $\Im(S,B)$ is HUR-stable on
$(S,B)$ with controls $(\psi(x), \frac{\psi(x)}{(1-L)|p(x)|})$ and
HUR-stable function $T$ is unique. The proof is complete.\\

With the Theorem of Z. Gajda (\ref{Gajda}), its easy to show that
the following result.
\begin{cor}\label{gaj}
Let $S=B$ be real field $\mathbb{R}$ and $\Im(S,B)$ be the
functional equation $f(2x)=2f(x)$. Then for every $\theta>0$ there
is no constant $\delta\in [0,\infty)$ such that the functional
equation $\Im(S,B)$ is HUR-stable on $(S,B)$ with controls
$(\theta(|x|),\delta|x|)$.
\end{cor}
Its obvious that the above corollary is a counterexample for the
special case of functional equation in the Theorem (\ref{6.1}),
when $L = 1$. With the Theorem (\ref{6.1}), we have the following
Corollary.
\begin{cor}
Let $\Im(S,B)$ be functional equation $f(\rho(x))=p(x)f(x)+q(x)$.
If $a\leq p(x)$ for all $x\in S$ and some real $a>1$, then for
every real $\delta>0$, the functional equation $\Im(S,B)$ is
HU-stable on $(X,B)$ with controls $(\delta,
\frac{a\delta}{a-1})$ and HU-stable function such $T$ is unique.
\end{cor}
\textbf{\emph{Proof.}} Sets $\psi(x):=\delta$ for all $x\in S$ and
$L:=\frac{1}{a}$. Now apply Theorem (\ref{6.1}).\\

Similarly, we prove that a Hyers-Ulam-Rassias stability for the
linear functional equation with another suitable conditions.
\begin{thm}\label{6.5}
Let $\Im(S,B)$ be functional equation $f(\rho(x))=p(x)f(x)+q(x)$.
Let there exists a positive real $L<1$ such that
\begin{eqnarray}
  |p(x)|\psi(\rho^{-1}(x))&\leq& L\psi(x)
\end{eqnarray}
for all $x\in S$ and also $\rho$ be a permutation of $S$. Then
the functional equation $\Im(S,B)$ is HUR-stable on $(S,B)$ with
controls $(\psi(x), \frac{1}{1-L}\psi(\rho^{-1}(x)))$ and
HUR-stable function such $T$ is unique.
\end{thm}
\emph{\textbf{Proof.}} Let us consider the set
$\mathcal{A}:=\{h:S\rightarrow B\}$ and introduce the generalized
metric on $\mathcal{A}$:
$$d(u,h)=\sup_{\{x\in S\ ;\ \psi(x)\neq0\}}\frac{\|g(x)-h(x)\|}{\psi(\rho^{-1}(x))}.$$
It is easy to show that $(\mathcal{A}, d)$ is generalized complete
metric space. Now we define the function $J
:\mathcal{A}\rightarrow \mathcal{A}$ with
$$J(h(x))=p(\rho^{-1}(x))h(\rho^{-1}(x))+q(\rho^{-1}(x))$$ for all $h\in \mathcal{A}$ and $x\in X$. Since
$|p(x)|\psi(\rho^{-1}(x))\leq L\psi(x)$ for all $x\in S$ and
$\rho$ is a permutation of $S$, so
  \begin{eqnarray*}
    d(J(u),J(h)) &=& \sup_{\{x\in S\ ;\ \psi(x)\neq0\}}\frac{|p(\rho^{-1}(x))|\|u(\rho^{-1}(x))-h(\rho^{-1}(x))\|}{\psi(\rho^{-1}(x))} \\
      &\leq& \sup_{\{x\in S\ ;\ \psi(\rho^{-1}(x))\neq0\}}L\frac{\|u(\rho^{-1}(x))-h(\rho^{-1}(x))\|}{\psi(\rho^{-2}(x))}\leq Ld(u,h)
  \end{eqnarray*}
for all $u, h\in \mathcal{A}$, that is $J$ is a strictly
contractive selfmapping of $\mathcal{A}$, with the Lipschitz
constant $L$ (note that $0<L<1$). We have
$$\|f(\rho(x))-p(x)f(x)-q(x)\|\leq \psi(x)$$
for all $x\in S$, we get
$$\|f(x)-p(\rho^{-1}(x))f(\rho^{-1}(x))+q(\rho^{-1}(x))\|\leq\psi(\rho^{-1}(x))$$
for all $x\in S$, which says that $d(J(f),f)\leq 1<\infty$. So, by
Theorem (\ref{3.1}), there exists a mapping $T:X\rightarrow B$
such that
  \begin{enumerate}
  \item $T$ is a fixed point of $J$, i.e.,
  \begin{eqnarray}\label{6.7}
    T(\rho(x)) &=& p(x)T(x)+q(x)
  \end{eqnarray}
  for all $x\in S$. The mapping $T$ is a unique fixed point of $J$ in the set $\tilde{\mathcal{A}}=\{h\in \mathcal{A}\ :\ d(f,h)<\infty\}$.
  This implies that $T\in Z_{\Im(S,B)}$ and is unique HUR-stable
  function. Also there exists $C\in (0,\infty)$ satisfying
  $$\|f(x)-T(x)\|\leq C\psi(\rho^{-1}(x))$$
  for all $x\in S$.
  \item $d(J^{n}(f),T)\rightarrow 0$ as $n\rightarrow \infty$. This implies that
  $$T(x)=\lim_{n\rightarrow\infty}  \prod_{i=1}^{n}p(\rho^{-i}(x))f(\rho^{-n}(x)) -\sum_{k=1}^{n} q(\rho^{i}(x)) \prod_{i=0}^{k-1}p(\rho^{-i}(x))$$
  for all $x\in S$.
  \item $d(f,T)\leq\frac{1}{1-L}d(J(f),f)$, which implies,
  $$d(f,T)\leq\frac{1}{1-L}.$$
  $$\|f(x)-T(x)\|\leq \frac{1}{1-L}\psi(\rho^{-1}(x))$$
for all $x\in S$.
\end{enumerate}
Therefore, the functional equation $\Im(S,B)$ is HUR-stable on
$(S,B)$ with controls $(\psi(x),
\frac{1}{1-L}\psi(\rho^{-1}(x)))$ and HUR-stable function $T$ is
unique. The proof is complete.\\

Similar to the Corollary (\ref{gaj}, we get the following result,
where its counterexample for the special case of functional
equation in the Theorem (\ref{6.5}), when $L = 1$.
\begin{cor}
Let $S=B$ be real field $\mathbb{R}$ and $\Im(S,B)$ be the
functional equation $f(\frac{1}{2}x)=\frac{1}{2}f(x)$. Then for
every $\theta>0$ there is no constant $\delta\in [0,\infty)$ such
that the functional equation $\Im(S,B)$ is HUR-stable on $(S,B)$
with controls $(\theta(|x|),\delta|x|)$.
\end{cor}
\begin{cor}
Let $\Im(S,B)$ be functional equation $f(\rho(x))=p(x)f(x)+q(x)$.
If $|p(x)|\leq L$ for all $x\in S$, some real $0<L<1$, and $\rho$
be a permutation of $S$. Then for every real $\delta>0$, the
functional equation $\Im(S,B)$ is HU-stable on $(S,B)$ with
controls $(\delta, \frac{\delta}{1-L})$ and HU-stable function
such $T$ is unique.
\end{cor}
\begin{cor}
Let $X$ be a normed linear space over $F$, let $\Im(S,B)$ be
functional equation $f(ax)=kf(x)$ for fixed constants $a$ and
$k$, and let $p\in \mathbb{R}$. If $p\leq0$, $|a|>1$ and $|k|>1$
or $p\leq0$, $|a|<1$ and $|k|<1$ or $p\geq0$, $|a|>1$ and $|k|<1$
or $p\geq0$, $|a|<1$ and $|k|>1$, then the functional equation
$\Im(S,B)$ is HUR-stable on $(S,B)$ with controls $(\|x\|^{p},
\frac{\|x\|^{p}}{||k|-1|})$ and HUR-stable function such $T$ is
unique.
\end{cor}
\textbf{\emph{proof.}} Set $\rho(x):=ax$ and $\psi(x):=\|x\|^{p}$
for all $x\in S$ and then apply Theorem (\ref{6.1}) and Theorem
(\ref{6.5}).\\

Now in the following we consider the Hyers-Ulam-Rassias stability
Pexiderized linear functional equation (\ref{e4}).
\begin{thm}\label{6.10}
Let $f,g:S\rightarrow B$ be a function and
\begin{eqnarray}\label{6.11}
  \|f(\rho(x))-p(x)g(x)-q(x)\| &\leq& \psi(x)
\end{eqnarray}
for all $x\in S$. If there exists a positive real $L < 1$ such
that
\begin{eqnarray}
  \psi(\rho(x)) &\leq& L|p(\rho(x))|\psi(x)\label{6.12};\\
  \|f(\rho(x))-g(\rho(x))\|&\leq&L\|f(x)-g(x)\| \label{6.13}
\end{eqnarray}
for all $x\in S$. Then there is an function $T$ such that
$T(\rho(x))= p(x)T(x)+q(x)$
$$\|f(x)-T(x)\|\leq \frac{\widetilde{\psi}(x)}{(1-L)|p(x)|}$$
$$\|g(x)-T(x)\|\leq \frac{L}{1-L}[\frac{\widetilde{\psi}(x)+\psi(x)}{|p(x)|}]$$
for all $x\in S$, in which
$\widetilde{\psi}(x)=\psi(x)+|p(x)|\|f(x)-g(x)\|$ for all $x\in
X$.
\end{thm}
\textbf{\emph{Proof.}} Applying (\ref{6.11}), we get
\begin{eqnarray*}
  \|f(\rho(x))-p(x)f(x)-q(x)\|&\leq&  \psi(x)+|p(x)|\|f(x)-g(x)\|\\
                             &\leq&\widetilde{\psi}(x)
\end{eqnarray*}
for all $x\in S$. From (\ref{6.12}) and (\ref{6.13}), its easy to
show that the following inequality
$$\widetilde{\psi}((\rho(x))\leq
L|p(\rho(x))|\widetilde{\psi}(x)$$ for all $x\in S$. So, by
Theorem (\ref{6.1}), there is an unique function $T: X\rightarrow
B$ such that $T(\rho(x))=p(x)T(x)+q(x)$
$$\|f(x)-T(x)\|\leq \frac{\widetilde{\psi}(x)}{(1-L)|p(x)|}$$
for all $x\in S$. So from last inequality, we have
$$\|f(\rho(x))-T(\rho(x))\|\leq  \frac{\widetilde{\psi}(\rho(x))}{(1-L)|p(\rho(x))|}$$
for all $x\in S$. We show that $T$ is a linear equation, thus from
last inequality and (\ref{6.11}), we get
$$\|g(x)-T(x)\|\leq \frac{L}{1-L}[\frac{\widetilde{\psi}(x)+\psi(x)}{|p(x)|}]$$
for all $x\in S$. The proof is complete.
\subsection{Common Stability for the Systems of Homogeneous Linear Equations}

Throughout this section, assume that $\{p_{i}:S\rightarrow
K\backslash\{0\}\}_{i\in I}$, $\{\rho_{i}: S\rightarrow S\}_{i\in
I}$ be two family of functions. Here $i$ is a variable ranging
over the arbitrary index set $I$. Also we define the functions
$P_{i,n}:S\rightarrow K\backslash\{0\}$ and
$\theta_{i,n}(x):X\rightarrow \mathbb{R}^{+}$ with
$$P_{i,n}(x)=\prod_{k=0}^{n-1}p_{i}(\rho_{i}^{k}(x))$$ and
$$\theta_{i,n}(x)=\frac{(1-L_{i}^{n})\psi_{i}(x)}{(1-L_{i})|p_{i}(x)|}$$
for a family of positive real's $\{L_{i}\}_{i\in I}$, all $x\in
S$, any index $i$ and positive integer $n$.

\begin{defn}\label{0hu-stable}
Let $S$ and $B$ be nonempty sets. Let $\{\psi_{i}:X\rightarrow
\mathbb{R}^{+}\}_{i\in I}$ be a family of functions with index
set $I$ , $\phi:S\rightarrow [0,\infty)$ be a function, and
$\{\Im_{i}(S,B)\}_{i\in I}$ be a family of functional equations.
If for every functions $f:S\rightarrow B$ satisfying the
inequality
\begin{equation}\label{00.2}
  d(G_{i}(f,x),G_{i}(f,x))\leq \psi_{i}(x)
\end{equation}
for all $x, y\in S$ and any $i\in I$, there exists $T\in
\bigcap_{i\in I} Z_{\Im_{i}(S,B)}$ such that
\begin{equation}\label{00.3}
  d(f(x),T(x))\leq \phi(x)
\end{equation}
for all $x\in S$, then we say that the family of functional
equations $\{\Im_{i}(S,B)\}_{i\in I}$ are common
Hyers-Ulam-Rassias stable on $(S,B)$ with control functions
$(\{\psi_{i}\}_{i\in I}, \phi)$ and we denoted it by
"CHUR-stable" on $(S,B)$ with controls $(\{\psi_{i}\}_{i\in I},
\phi)$. Also we call the function $T$ as "CHUR-stable function".
\end{defn}
In this section, we consider some systems of homogeneous linear
equations
\begin{equation}\label{6.14}
    f(\rho_{i}(x))=p_{i}(x)f(x),
\end{equation}
and our aim is to establish some common Hyers-Ulam-Rassias
stability for these systems of functional equations. As a
consequence of these results, we give some applications to the
study of the superstability result for exponential functional
equation to the a family of functional equations. Note that the
following Theorem is partial affirmative answer to problem 1, in
the 13st ICFEI.
\begin{thm} \label{6.15}
Let $\{\Im_{i}(S,B)\}_{i\in I}$ be a family of functional
equations, in which $\Im_{i}(S,B)$ be the functional equation
$f(\rho_{i}(x))=p_{i}(x)f(x)$ for any $i\in I$. Assume that
\begin{enumerate} \label{6.17}
    \item there exists a family of positive real's $\{L_{i}\}_{i\in I}$ such that
    $L_{i}<1$ and
    $$\psi_{i}(\rho_{i}(x))\leq
    L_{i}|p_{i}(\rho_{i}(x))|\psi_{i}(x)$$ for all $x\in S$ and $i\in
    I$;
    \item $\rho_{i}\rho_{j}=\rho_{i}\rho_{j}$ for all $i,j\in I$;
    \item $p_{i}(\rho_{j}(x))=p_{i}(x)$ for all distinct $i,j\in I$,
    ;
    \item $\lim_{n\rightarrow\infty} \frac{\theta_{i,n}(\rho_{j}^{n}(x))}{|P_{j,n}(x)|}=0$ for all
    $x\in X$ and every distinct $i,j\in I$.
\end{enumerate}
Then the family of functional equations $\{\Im_{i}(S,B)\}_{i\in
I}$ "CHUR-stable" on $(S,B)$ with controls $(\{\psi_{i}\}_{i\in
I}, \inf_{i\in I}\{\frac{\psi_{i}(x)}{(1-L_{i})|p_{i}(x)|}\})$ and
the CHUR-stable function such $T$ is unique.
\end{thm}
\textbf{\emph{Proof.}} According to the our auumptions, for every
$i\in I$, the Theorem (\ref{6.1}), implies that the functional
equation $\Im_{i}(S,B)$ is HUR-stable on $(S,B)$ with controls
$(\psi(x), \frac{\psi(x)}{(1-L)|p(x)|})$ and HUR-stable function
such $T_{i}\in Z_{\Im_{i}(S,B)}$ is unique. Moreover, The
function $T_{i}$ is given by
$$T_{i}(x)=\lim_{n\rightarrow\infty}\frac{f(\rho_{i}^{n}(x))}{\prod_{k=0}^{n-1}p_{i}(\rho_{i}^{k}(x))}=\lim_{n\rightarrow\infty}J_{i}^{n}(f)$$
for all $x\in S$ and any fixed $i\in I$. In the proof of Theorem
(\ref{6.1}), we show that
$$d(J_{i}(f),f)\leq 1.$$
By induction, its easy to show that
$$d(J_{i}^{n}(f),f)\leq\frac{1-L_{i}^{n}}{1-L_{i}},$$
which says that
$$\|f(\rho_{i}^{n}(x))-\prod_{k=0}^{n-1}p_{i}(\rho_{i}^{k}(x))f(x)\|\leq
(\prod_{k=0}^{n-1}|p_{i}(\rho_{i}^{k}(x))|)\frac{(1-L_{i}^{n})\psi_{i}(x)}{(1-L_{i})|p_{i}(x)|}$$
for all $x\in S$ and $i\in I$. Now we show that $T_{i}=T_{j}$ for
any $i,j\in I$. Let $i$ and $j$ be two arbitrary fixed indexes of
$I$. So, from last inequality, we obtain
\begin{eqnarray}
  \|f(\rho_{i}^{n}(x))-P_{i,n}(x)f(x)\| &\leq& |P_{i,n}(x)|\theta_{i,n}(x) \label{6.18} ;\\
  \|f(\rho_{j}^{n}(x))-P_{j,n}(x)f(x)\| &\leq& |P_{j,n}(x)|\theta_{j,n}(x)
  \label{6.19}
\end{eqnarray}
for all $x\in S$. On the replacing $x$ by $\rho_{j}^{n}(x)$ in
(\ref{6.18}) and $x$ by $\rho_{i}^{n}(x)$ in (\ref{6.19})
\begin{eqnarray*}
  \|f(\rho_{i}^{n}(\rho_{j}^{n}(x)))-P_{i,n}(\rho_{j}^{n}(x))f(\rho_{j}^{n}(x))\| &\leq& |P_{i,n}(\rho_{j}^{n}(x))|\theta_{i,n}(\rho_{j}^{n}(x));\\
  \|f(\rho_{j}^{n}(\rho_{i}^{n}(x)))-P_{j,n}(\rho_{i}^{n}(x))f(\rho_{i}^{n}(x))\| &\leq& |P_{j,n}(\rho_{i}^{n}(x))|\theta_{j,n}(\rho_{i}^{n}(x))
\end{eqnarray*}
for all $x\in X$. From assumptions (2) and (3), its obvious that
$f(\rho_{i}^{n}(\rho_{j}^{n}(x)))=f(\rho_{j}^{n}(\rho_{i}^{n}(x)))$
, $P_{i,n}(\rho_{j}^{n}(x))=P_{i,n}(x)$ and
$P_{j,n}(\rho_{i}^{n}(x))=P_{j,n}(x)$ for all $x\in S$. So, from
last two inequalities
$$\|P_{i,n}(x)f(\rho_{j}^{n}(x))-P_{j,n}(x)f(\rho_{i}^{n}(x))\|\leq |P_{i,n}(x)|\theta_{i,n}(\rho_{j}^{n}(x))+|P_{j,n}(x)|\theta_{j,n}(\rho_{i}^{n}(x))$$
or
$$\|\frac{f(\rho_{j}^{n}(x))}{P_{j,n}(x)}-\frac{f(\rho_{i}^{n}(x))}{P_{i,n}(x)}\|\leq \frac{\theta_{i,n}(\rho_{j}^{n}(x))}{|P_{j,n}(x)|}+\frac{\theta_{j,n}(\rho_{i}^{n}(x))}{|P_{i,n}(x)|}$$
for all $x\in X$. From assumption $\lim_{n\rightarrow\infty}
\frac{\theta_{i,n}(\rho_{j}^{n}(x))}{|P_{j,n}(x)|}=0$ for all
$x\in S$ and every distinct $i,j\in I$, so, its implies that
$T_{i}=T_{j}$.

Now set $T=T_{i}$ and since $\|f(x)-T_{i}(x)\|\leq
\frac{\psi_{i}(x)}{(1-L_{i})|p_{i}(x)|}$ for all $x\in S$ and all
$i\in I$, there is a unique function $T$ such that
$$T(\rho_{i}(x))=p_{i}(x)T(x)$$ (i.e., $T\in
\bigcap_{i\in I} Z_{\Im_{i}(S,B)}$) for all $x\in S$ and $i\in I$
and also
 $$\|f(x)-T(x)\|\leq \inf_{i\in
    I}\{\frac{\psi_{i}(x)}{(1-L_{i})|p_{i}(x)|}\}$$
for $x\in S$. The proof is complete.
\begin{cor}\label{6.20}
Let $\{\Im_{i}(S,B)\}_{i\in I}$ be a family of functional
equations, in which $\Im_{i}(S,B)$ be the functional equation
$f(\rho_{i}(x))=c_{i}f(x)$, where $\{c_{i}\}_{i\in I}$ is a
family of constants. Assume that
$\rho_{i}\rho_{j}=\rho_{i}\rho_{j}$ for all $i,j\in J$ and also
there exists a family of positive real's $\{L_{i}\}_{i\in J}$
such that $0<L_{i}<1$ and
\begin{equation}\label{qwe}
  \psi_{i}(\rho_{i}(x))\leq
    L_{i}|c_{i}|\psi_{i}(x)
\end{equation}
for all $x\in S$ and $i\in J$, in which $J=\{i\in I :\
|c_{i}|>1,\ L_{i}|c_{i}|\in(0,1]\}$, then the family of
functional equations $\{\Im_{i}(S,B)\}_{i\in J}$ "CHUR-stable" on
$(S,B)$ with controls $(\{\psi_{i}\}_{i\in J}, \inf_{i\in
J}\{\frac{\psi_{i}(x)}{(1-L_{i})c_{i}}\})$ and the CHUR-stable
function such $T$ is unique.
\end{cor}
\textbf{\emph{Proof.}} Set $p_{i}(x):= c_{i}$ for all $i\in J$,
then the conditions (1), (2), and (3) of Theorem (\ref{6.15}) is
holds and now if we show that the condition (4) is hold, then
with we obtain the result. We have
$$\theta_{i,n}(\rho_{j}^{n}(x))=\frac{(1-L_{i}^{n})\psi_{i}(\rho_{j}^{n}(x))}{(1-L_{i})|c_{i}|}$$
and
$$P_{j,n}(x)=c_{j}^{n}$$
for all $x\in S$ and and every $i, j\in J$. Now from (\ref{qwe}),
we get $\psi_{i}(\rho_{j}^{n}(x))\leq
L_{i}^{n}|c_{i}|^{n}\psi_{i}(x)$ and so that
$$\lim_{n\rightarrow\infty}
\frac{\theta_{i,n}(\rho_{j}^{n}(x))}{|P_{j,n}(x)|}\leq\frac{(1-L_{i}^{n})L_{i}^{n}|c_{i}|^{n}\psi_{i}(x)}{|c_{j}|^{n}}$$
for any $x\in S$ and every distinct $i,j\in J$. Since
$0<|L_{i}c_{i}|<1$ and $|c_{i}|>1$ for all $i\in J$, so the above
limit approach to zero. The proof is complete.\\

Now with the above Corolary, in the following, we prove a
superstability result for Cauchy exponential functional equation,
where we discussed about it in 2th section
\begin{thm}
Let $S$ be commutative semigroup and $B$ be the field of complex
numbers $\mathbb{C}$. Let $\Im(S,B)$ be the functional equation
$f(x\cdot y)=g(y)f(x)$, where $g:S\rightarrow C$ is a function.
Let the set $J$ be the elements of $i\in S$, where $|g(a)|>1$ and
and there exists $L_{i}\in (0,1)$ with $L_{i}|g(i)|\in(0,1]$.
Assume that $\phi:S^{2}\rightarrow \mathbb{R}^{+}$ is function,
$g$ be is unbounded function
$$\phi(x,y\cdot i)\leq L_{a}|g(a)|\phi(x,y)$$
for all $x, y\in S$ and $i\in J$. If $g(J)$ be an unbounded set,
then the functional equation $\Im(S,B)$ is hyperstable on $(S,B)$
with control $(\phi)$.
\end{thm}
\textbf{Proof.} Let $g$ be a unbounded function, then sets
$\rho_{i}(x):=x\cdot i$, $c_{i}:=g(i)$, and $\psi_{i}:=\phi(i,x)$
for all $x\in S$ and any $i\in J$. Since
$\rho_{i}\rho_{j}=\rho_{i}\rho_{j}$ and
$\psi_{i}(\rho_{i}(x))\leq L_{a}|g(a)|\psi_{i}(x)$ for all $x\in
S$ and any $i\in J$, then by Corollary (\ref{6.20}), the family
of functional equations $\{\Im_{i}(f)\}_{i\in J}$ "CHUR-stable"
on $(S,B)$ with controls $(\{\psi_{i}\}_{i\in J}, \inf_{i\in
J}\{\frac{\psi_{i}(x)}{(1-L_{i})c_{i}}\})$ and the CHUR-stable
function such $T$ is unique. Since $g$ is a unbounded function,
from last inequality $T=f$ (note that $L_{i}|g(i)|\in(0,1]$),
which implies that
$$f(\rho_{i}(x))=c_{i}f(x)$$
or
\begin{equation}\label{6.22}
    f(x\cdot i)=g(i)f(x)
\end{equation}
for all $x\in S$ and $i\in J$. We have
\begin{equation}\label{gh1}
 \|f(x\cdot y)-g(y)f(x)\|\leq\phi(x,y)
\end{equation}
for all $x,y \in S$. On the replacing $y$ by $y\cdot i^{n}$ in
(\ref{gh1})
$$\|f((x\cdot y)\cdot i^{n})-g(y\cdot
i^{n})f(x)\|\leq\phi(x,y\cdot i^{n})$$ or
$$\|\frac{f((x\cdot y)\cdot i^{n})}{g(i)^{n}}-\frac{g(y)f(x\cdot
i^{n})}{g(i)^{n}}\|\leq\frac{\phi(x,y\cdot i^{n})}{|g(i)|^{n}}$$
for all $x,y\in S$, any fixed $i\in I$ and positive integer $n$.
From equation (\ref{6.22}), its easy to show that $f(x\cdot
i^{n})=g(i)^{n}f(x)$ and $\phi(x,y\cdot i^{n})\leq
L_{i}^{n}|g(i)|^{n}\phi(x,y)$ for all $x\in S$, any fixed $i\in J$
and positive integer $n$. So, we have
$$\|f(x\cdot y)-g(y)f(x)\|\leq L_{i}^{n}\phi(x,y)$$
for all $x,y\in S$, any fixed $i\in J$ and positive integer $n$
), which implies that $f(x\cdot y)=g(y)f(x)$ for all $x,y\in S$
(i.e., $f\in Z_{\Im(S,B)}$). The proof is complete.
\section{\textbf{Acknowledgments}}
When, I was undergraduate student of pure mathematics at
university of Mazandaran (2008-2013), I was interesting to
research in mathematics with problem solving specially about
"Functional Equation" and then Professor Mohsen Alimohammady
helped me for research in mathematics. I'm very much grateful and
dedicate this thesis to Professor.


\bibliographystyle{amsplain}


\bigskip
\bigskip

{\footnotesize \pn{\bf Ali Sadeghi}\; \\ {Graduate Student of
Pure Mathematics at Department of Mathematics}, {Tarbiat Modares University, P.O.Box 14115-333,} {Tehran, Iran}\\
{\tt Email: sadeghi.ali68@gmail.com; ali.sadeghi@modares.ac.ir}\\

\end{document}